\documentclass[11pt]{article}

\usepackage{color}
\usepackage{amsmath}
\usepackage{mathrsfs}
\usepackage{amssymb}
\usepackage{amsthm}
\usepackage{epstopdf}
\usepackage{graphicx}
\usepackage{mathtools}
\usepackage[usenames,dvipsnames, table]{xcolor}
\usepackage[hang]{caption}
\usepackage{hyperref}
\usepackage{curve2e}

\usepackage{oldgerm}  
\usepackage[yyyymmdd,hhmmss]{datetime}

\usepackage{mathabx} 
\def\bee{\begin{enumerate}}\def\eee{\end{enumerate}}
\def\bei{\begin{itemize}}\def\eei{\end{itemize}}
\newcommand{\nco}{\newcommand}
\def\R{\mathbb{R}}\def\C{\mathbb{C}}
\nco{\red}{\color{red}}
\nco{\blue}{\color{blue}}
\nco{\green}{\color{green}}
\nco{\cyan}{\color{cyan}}
\nco{\brown}{\color{Magenta}}
\def\Blue#1{\blue #1\normalcolor}

\nco{\magenta}{\color{magenta}}

\nco{\violet}{\color{violet}}
\nco{\olive}{\color{Emerald}}
\nco{\orange}{\color{orange}}
\nco{\redend}{\normalcolor}
\nco{\blueend}{\normalcolor}
\def\inv#1{\frac{1}{#1}}
\def\tr{{\rm tr}\,}
\def\sign{{\rm sign}\,}

\def\ommit#1{{}}
\def\({\left(}
\def\){\right)}
\def\ie{{\it i.e.,\/}\ }
\def\ie{{\rm i.e.,\/}\ }

\definecolor{cb}{rgb}{.8,.5,0}

\nco{\rnc}{\renewcommand}
\rnc{\title}[1]{{\Large\bf\mbox{}\\\medskip#1\bigskip\medskip\\}}
\rnc{\author}[1]{{\large #1\smallskip\\}}
\nco{\address}[1]{{\em #1\medskip\\}}

\def\diag{{\rm diag \,}}
\def\ii{\mathrm{i\,}}
\nco{\bun}{{\bf 1}}
\def\be{\begin{equation}}\def\ee{\end{equation}}
\def\bea{\begin{eqnarray}}\def\eea{\end{eqnarray}}
\def\bee{\begin{enumerate}}\def\eee{\end{enumerate}}
\def\bei{\begin{itemize}}\def\eei{\end{itemize}}
\def\oh{\frac{1}{2}}
\def\ommit#1{{}}

\def\SU{{\rm SU}}\def\SO{{\rm SO}}
\def\inv#1{\frac{1}{#1}}
\def\mult{{\rm mult}}
\def\tr{{\rm tr\, }}

\def\eq=#1{\buildrel #1 \over{=}}
\def\KT{KT\ }
\def\GT{GT\ }

\def\CH{{\mathcal H}}  \def\CI{{\mathcal I}}   \def\CJ{{\mathcal J}}   \def\CO{{\mathcal O}}
\def\CP{{\mathcal P}}

\def\Gb{\pmb{\beta}}

\def\p{\mathrm{p}}

\def\diag{{\rm diag \,}}
\def\ii{\mathrm{i\,}}
\def\N{\mathbb{N}}

\def\R{\mathbb{R}}
\def\C{\mathbb{C}}

\begin{document}
\begin{titlepage}
%
\begin{center}
\title{On Schur problem\\[6pt]
{and Kostka numbers}}
\medskip
\author{Robert Coquereaux} 
\address{Aix Marseille Univ, Universit\'e de Toulon, CNRS, CPT, Marseille, France\\
Centre de Physique Th\'eorique}
\bigskip\medskip
\medskip
\author{Jean-Bernard Zuber}
\address{
\ Sorbonne Universit\'e,  UMR 7589, LPTHE, F-75005,  Paris, France\\ \& CNRS, UMR 7589, LPTHE, F-75005, Paris, France
 }

%
\bigskip\bigskip
\begin{abstract}
 We reconsider the two related problems: distribution of the diagonal elements of a Hermitian $n\times n$ 
matrix of known eigenvalues (Schur)
and determination of multiplicities of weights in a given irreducible representation of $\SU(n)$ (Kostka). It is well known that 
the former yields a semi-classical picture
of the latter. We present explicit expressions for low values of $n$ that complement those
given in the literature \cite{GLS,BGR}, recall some exact (non asymptotic) relation between the two problems,  comment on the 
limiting procedure whereby Kostka numbers are obtained from Littlewood--Richardson coefficients, and finally extend these
considerations to the case of the $B_2$ algebra, with a few novel conjectures.
\end{abstract}
\end{center}

 \end{titlepage}

 \section{Introduction}
 \subsection{Overview}
It has been known for long that multiplicity problems in representation theory have an asymptotic limit that may be treated by 
semi-classical methods \cite{He82, GLS}. For example, the behavior for large representations of generalized Littlewood--Richardson 
(LR) coefficients, \ie coefficients of decomposition into irreducible representations (irreps) of the tensor product of two irreps, 
admits a semi-classical description in terms of Horn's problem.
For a review and a list of references, see \cite{CMSZ1}. 
Similarly, one may consider the Kostka 
numbers $\mult_\lambda(\delta)$, \ie the multiplicity of weight $\delta$ in the irrep of highest weight $\lambda$.
As also well known \cite{GLS,BGR},  the asymptotics (for large $\lambda$ and $\delta$) of those numbers 
is related to Schur's problem, which deals with the properties and distribution of diagonal elements of a Hermitian matrix of given spectrum.

 In the following, we first reexamine that classical Schur problem: 
  in the case where the original Hermitian matrix is taken at random 
uniformly on its orbit, we recall (sect. \ref{sec2}) how the probability density function (PDF) of its diagonal elements is determined by a 
{\it volume function $\CI$}, which has as a support the permutahedron determined by the eigenvalues  and is a piecewise polynomial 
function of these eigenvalues and of the diagonal elements, with non analyticities of a prescribed type on 
an a priori known locus.  We give quite explicit expressions of that function for the cases of
$\SU(3)$ and $\SU(4)$ (coadjoint) orbits (sect. \ref{lowvalues}). 
In sect. \ref{sec4}, we turn to the parallel representation-theoretic problem, namely the determination 
of Kostka numbers. We rederive (sect. \ref{CI-mult}) in this new context an exact  relationship between these multiplicities and the aforementioned 
volume function, which was already discussed in the Horn problem \cite{CZ1,CMSZ1}. That relation leads in a natural way to the
semi-classical asymptotic   limit mentioned above  (sect. \ref{Asymptotics}). 
A guiding thread through this work is the (fairly obvious and well known) fact that Kostka numbers may be obtained as a 
certain limit of Littlewood--Richardson coefficients when two of their arguments grow large with a finite difference. In fact that 
limit is approached quite fast, and it is an intriguing problem to find values for the threshold value $s_c$ of the scaling parameter 
beyond which the asymptotic Kostka number is reached. 
We address that question in sect. \ref{KostkaFromLR} and propose a  (conservative) upper bound on that threshold in the case of SU(4).
We also recall the combinatorial interpretation 
of Kostka's numbers, in terms of reduced Knutson--Tao honeycombs aka Gelfand--Tsetlin triangles, or of reduced O-blades (sect. \ref{polytope}, \ref{ObladesVersusLianas}). 
The latter may be recast in a new picture of ``forests of lianas", as discussed in the Appendix. 
Finally, sect. \ref{B2} is devoted to a discussion of what can be said or conjectured in the case of the $B_2$ case.

Several of these results have already appeared in some guise in the literature. 
The domains of polynomiality of the Duistermaat--Heckman  measure  and the transitions between them have been discussed for 
$\SU(n)$, with emphasis on  $n=3$ and 4, in \cite{GLS}. The parallel analysis of multiplicities has been carried out in \cite{BGR}, 
using the method of vector partition functions. A general discussion making use of Littelmann's paths \cite{Li} has been done by Bliem \cite{Bl},
with an illustration in the case of the $B_2=\mathfrak{so}(5)$ algebra.
We believe, however, that several aspects of our approach are
original, that our results for the case $\SU(4)$, resp. $B_2$,  complement those of \cite{GLS, BGR}, resp. \cite{Bl}, 
that the discussion of the threshold value $s_c$ is novel, 
and that the liana picture may give a new insight in the combinatorial aspects of the problem. 

This Schur problem --as its sibling the Horn problem-- presents a unique and fascinating mix  of various ingredients, 
algebraic, geometric and group theoretic. As such, we hope that our modest contribution would have pleased our distinguished 
colleague Boris Dubrovin, who, all his life, paid an acute attention to the interface 
between mathematics and physics.


\subsection{Notations}
\label{sec:notations}
In the following, we use two alternative notations for the objects pertaining to $\SU(n)$.
First in the classical Schur problem,  the ordered eigenvalues of a $n\times n$ 
Hermitian matrix will be denoted by $\alpha=(\alpha_1,\cdots, \alpha_n)$, with round brackets, 
and its diagonal elements by $\xi=(\xi_1,\cdots, \xi_n)$, with clearly $\sum_i \alpha_i=\sum_i \xi_i$.  
 At the possible price of an overall shift of all $\alpha_i$, one may assume they are all non negative.
In the case they are all non negative integers,  one may regard them as defining a partition and encode them in a Young diagram.

A vector $\lambda$ of the weight lattice of $\SU(n)$ may be denoted 
by its $(n-1)$ (Dynkin) components in the fundamental weight basis: $\lambda=\{ \lambda_1,\cdots, \lambda_{n-1} \}$,
with curly brackets, or as a partition, with $n$ ``Young components"
 equal to the lengths of rows of the corresponding Young diagram:  $\alpha=(\alpha_1,\cdots, \alpha_n)$. The former
is recovered from the latter by $\lambda_i=\alpha_i-\alpha_{i+1},\ \ i=1,\cdots,n-1$. 

 Conversely, when dealing with the  highest weight  $\lambda$ of an irreducible representation (irrep) of $\SU(n)$,
it is natural to define a decreasing partition $\alpha$ with $\alpha_i=\sum_{j=i}^{n-1} \lambda_{j},\ \  i=1,\cdots,n-1$,
and $\alpha_n=0$. For a weight  $\delta$ of that irrep, one  defines the sequence $\xi = (\xi_1,\cdots, \xi_n)$,  also called weight, 
with $\xi_i=\sum_{j=i}^{n-1} \delta_j+c,\ \  i=1,\cdots,n-1$, $\xi_n=c$, and  
one chooses $c=\inv{n} \sum_{i=1}^{n-1} {i} (\lambda_i-\delta_i)$, an integer,  in such a way that $\sum_{i=1}^n \alpha_i=\sum_{i=1}^n \xi_i$. 
 So $\xi$ is a non necessarily decreasing partition of  $\sum_{i=1}^n \alpha_i$.
 It is well known that $\mult_\lambda(\delta)$ is equal to the number of $\SU(n)$  semi-standard Young tableaux with fixed shape $\alpha$ and weight $\xi$.
In what follows, we shall use both languages interchangeably, and use the notations $(\lambda, \delta)$ or $(\alpha, \xi)$ with the above meaning, without further ado, writing for instance $\mult_\lambda(\delta) = \mult_\alpha(\xi)$.
Actually, we shall use the same notations and conventions, even if $\lambda$ and $\delta$ are not integral (so that $\alpha$ and $\xi$ are no longer partitions in the usual sense).
We hope that the context will prevent possible confusions.


\section{Schur's problem}
  \label{sec2}
Schur's problem deals with the following question:  If  $A$ is a $n$-by-$n$ Hermitian matrix with known eigenvalues $\alpha_1 \ge \hdots \ge \alpha_n$, what can be said about the diagonal elements $\xi_i:= A_{ii} $, $i=1,\cdots,\,n$?   As shown by Horn \cite{Ho54}, 
the $\xi$'s  lie in the closure of the {\it permutahedron} $\CP_\alpha$, \ie, the  convex polytope in $\R^{n}$ whose vertices are 
the points $(\alpha_{P(1)}, \alpha_{P(2)},\cdots, \alpha_{P(n)})$, $P\in \mathfrak{S}(n)$. 
 Note that by a translation of $A$ by a multiple of the identity matrix, we could always manage  to have 
 $\sum_i \xi_i=\sum_i \alpha_i=\tr A=0$~\footnote{This would be quite natural since the Lie algebra of $\SU(n)$ is the set of traceless skew-Hermitian matrices.},
  but we shall not generally assume this tracelessness in the following.

A more specific question is the following: if $A$ is drawn randomly and uniformly on its $\SU(n)$-orbit $\CO_\alpha$, what is the PDF of the $\xi$'s ?
It turns out that this PDF is, up to a factor,  the (inverse) Fourier transform of the orbital integral, \ie
the density of Duistermaat--Heckman's measure.  To show that, we follow the same steps as in \cite{Z1}.  
 The characteristic function of the random variables $\xi$, \ie the  Fourier transform of the desired PDF, is 
 $$ \varphi(x)  =\Bbb{E} (e^{\ii \sum_j x_j A_{jj}}) =\int dU \exp \ii \sum_{j=1}^n x_j \, (U.\alpha.U^\dagger)_{jj} 
 \,$$
 with $dU$ the normalized $\SU(n)$ Haar measure, and $x$ belongs to the $ \R^{n-1}$ hyperplane defined as $\sum_j x_j=0$.
 From this $\varphi(x) $ we recover the PDF $\p$ of $\xi$ by inverse Fourier transform
 $$ \p(\xi|\alpha)= \int \frac{d^{n-1}x}{(2\pi)^{n-1}}\, dU\, \exp \ii \sum_j x_j ((U.\alpha.U^\dagger)_{jj} -\xi_j)$$
 which is indeed the (inverse) Fourier integral of the orbital (HCIZ) integral $\CH(\alpha,\ii x) =\int dU e^{\ii \tr (xU\alpha U^\dagger)}$
 \bea \label{invFtr1} \p(\xi|\alpha)  &=&  \int \frac{d^{n-1}x}{(2\pi)^{n-1}}\, e^{- \ii \sum_j x_j \xi_j}\, \CH(\alpha,\ii x) 
\\ \label{invFtr}  &=&  \frac{\Delta(\rho)}{\Delta(\alpha)} \int\frac{d^{n-1}x}{(2\pi)^{n-1}}\,\inv{\Delta(\ii x)} \sum_{w\in \mathfrak{S}(n)} \epsilon(w) e^{\ii \sum_j x_j (\alpha_{w(j)} -\xi_j)}\,. \eea
In the above expressions $\epsilon(w)$ is the signature of $w$,
 $\Delta$ stands for the Vandermonde determinant, ${\Delta(\alpha)} =\prod_{1\le i<j\le n} (\alpha_i-\alpha_j)$
and 
 $\rho$ denotes the Weyl vector, written here in partition components $\rho=(n-1,n-2,\dots, 0)$.  In the present case of ${su}(n)$, 
 $\Delta(\rho)$ is the superfactorial $\prod_{j=1}^{n-1} j!$. For a general simple Lie algebra, 
 $\Delta(\xi)=\prod_{\alpha>0}\langle \alpha,\xi\rangle$, 
 a product over the positive roots; in the simply-laced cases $\Delta(\rho)$ 
 is the product\footnote{In the non simply-laced cases one has to divide this product by  
 an appropriate scaling coefficient (see for instance \cite{CMSZ1}) which is equal to $2^r$ for $B_r$, so that for  $B_2$ considered in sect. \ref{B2},  one gets $\Delta(\rho)= (1! \times 3! ) (2^2) = 3/2$.} 
   of factorials of the Coxeter exponents of the algebra.

\smallskip
  {\small It turns out that this may be recovered heuristically in a different way.
  Recall the connection between Schur's and Horn's problems~\cite{CMSZ2}. Consider Horn's problem for two matrices $A\in \CO_\alpha$ and $B\in\CO_\beta$ and assume that $\alpha \ll \beta$ and that the $\beta_j$ are distinct, $\beta_j-\beta_{j+1}\gg 1$. The eigenvalues of $C=A+B$,  to the first order of perturbation
 theory, are of the form $\gamma_i=\beta_i + A_{ii}= \beta_i+\xi_i$. The $\beta$'s being given, the  PDF of the $\gamma$'s is the PDF 
 of the $\xi$'s. The former, namely
 \be\label{palt} \p(\gamma| \alpha,\beta) = \inv{n!}  \frac{\Delta(\gamma)^2}{\Delta(\rho)^2 }  \int \frac{d^{n-1}x}{(2\pi)^{n-1}} \, \Delta(x)^2\, \CH(\alpha,\ii x) \CH(\beta,\ii x) 
 \CH(\gamma, -\ii x) \ee
 reduces in the limit to (\ref{invFtr1}). 
 To prove it, we expand as usual $\p(\gamma| \alpha,\beta)$ as
 \be\label{paltp}\p(\gamma| \alpha,\beta) =  \frac{\Delta(\gamma)\Delta(\rho)}{\Delta(\alpha) \Delta(\beta)}  \int \frac{d^{n-1}x}{(2\pi)^{n-1}}\, \inv{\Delta(\ii x)}\,\sum_{w,w'\in \mathfrak{S}(n)}  \epsilon(w . w')  e^{\ii \sum_j x_j (\alpha_{w(j)} +\beta_{w'(j)}-\gamma_j)}\ee
  and notice that for $\alpha.x\sim O(1)$, $\beta .x \gg 1$, all terms $e^{\ii x. \beta_{w'}}, \ w'\ne 1$ are rapidly oscillating and are suppressed,
  leaving only the term $w'=1$ for which the exponential reduces to $e^{\ii  \sum_j x_j (\alpha_{w(j)} -\xi_j)}$ 
  while $\Delta(\gamma)/\Delta(\beta)\approx 1$, thus reducing (\ref{paltp}) to (\ref{invFtr}). }

   In these expressions, the integration is over the hyperplane $\R^{n-1}$ defined by $\sum_i x_i=0$ (in fact, the Cartan algebra of $\SU(n)$).
 As usual, we change variables $u_i=x_i-x_{i+1}$, and denote $\tilde\Delta(u):=\Delta(x)$. We thus write
\bea \label{p-CI}   \p(\xi|\alpha)&=&  \frac{\Delta(\rho)}{\Delta(\alpha)} \CI(\alpha;\xi)\\
\label{CI}
\CI(\alpha;\xi)&=& \sum_{w\in \mathfrak{S}(n)} \epsilon(w)\int \frac{d^{n-1}u}{\tilde\Delta(\ii u)}    e^{\ii \sum_{j=1}^n  u_j  \sum_{k=1}^j (\alpha_{w(k)}-\xi_k)} \,,\eea
 and focus our attention on that function $\CI$, that, for reasons explained later, we may call ``the volume function of the Schur--Kostka problem''.
As everywhere in this paper we freely use notations that may refer either to Dynkin components (weights: $\lambda, \delta$) or to Young components (partitions: $\alpha, \xi$), so that with the conventions defined in (sect. \ref{sec:notations}) we may write,  for instance,  $\CI(\lambda;\delta)=\CI(\alpha;\xi)$.

It is first clear that, by definition, $\p$ and $\CI$ must be symmetric functions of the $\xi_i$, $i=1,\cdots,n$.

 By {\it a priori} arguments \cite{He82, DH}, or by explicit computation of (\ref{CI}), it is clear that $\CI$ is a piece-wise polynomial of its 
 arguments $\alpha$ and $\xi$, homogeneous of degree $(n-1)(n-2)/2$. By Riemann--Lebesgue theorem (\ie looking at  the decay of
 the integrand of  (\ref{CI}) at large $u$), 
 its differentiability class is $C^{n-3}$, just like in the parallel Horn's problem. 
 
 By using the same arguments as in \cite{CMSZ1}, or by applying the previous limiting procedure to the $\CJ$ function of Horn's problem, 
 we can assert {\it a priori} that the loci of non-analyticity of $\CI$ (\ie the places where its polynomial determination changes) {are  contained in} 
 the hyperplanes
 \be \xi_i=\alpha_{w(i)},\ \mathrm{or} \ \xi_i+\xi_j=\alpha_{w(i)}+\alpha_{w(j)},\ 
 \qquad w\in \mathfrak{S}(n)\,, \ee
 or more generally 
 \be\label{equsing} \sum_{i\in I} \xi_i  =\sum_{j\in J} \alpha_j\ee
with $I,J\subset \{1,2,\cdots, n\}, \ |I|=|J|\le \lfloor \frac{n}{2}\rfloor$.

The relationship between the degree of non-differentiability and the cardinality of the sets $I$ and $J$ in 
 \ref{equsing}
has been addressed in the fundamental work of  
Guillemin--Lerman--Sternberg~\cite{GLS} \S 3.5 . It is proven there that the ``jump" (\ie the change of determination) 
of $\CI$ across a singular hyperplane is of the form
\be\label{jumpGLS} \Delta \CI=  \mathrm{const.} \ (\sum_I \xi_i -\sum_J \alpha_j)^{m-1}+\cdots \ee
{where $m$ is an integer depending on $k = |I| = |J|$, which we will now determine}, 
thus showing that $\CI$ is of class $C^{m-2}$ across that singularity. 
Let's compute that number $m-1$ in the case of a singular hyperplane of the form (\ref{equsing}).
A generic orbit of $\SU(n)$ has dimension $d_n:= (n^2-1) -(n-1)=n(n-1)$. Let $k=|I|=|J|$,  with $k \le \lfloor \frac{n}{2}\rfloor$,
then 
\be\label{myster} m-1  = (d_n -d_k-d_{n-k})/2-1 =n k -k^2 -1=k(n-k)-1\,.\ee
 Thus for $\SU(4)$, $k=1$, resp. 2, leads to $m-1=  2$, resp. 3,  \ie corresponds to a $C^1$, resp. $C^2$ differentiability class;
 for $\SU(5)$, $k=1$, resp. 2, leads to $m-1=  3$, resp. 5, and differentiability class $C^2$, resp. $C^4$, etc.
 This will be fully corroborated by the explicit expressions given below.

In parallel to the study 
  of the ``volume function" $\CI$, one may consider the  properties of the multiplicity function $\mult_\lambda(\delta)$.
  It is known that it is also a piece-wise polynomial function of $\lambda$ and $\delta)$ \cite{He82}. Quite remarkably, it has
  been proved that its singularities (changes of polynomial determination) as a function of $\delta$ occur on the same locus as
  those of $\CI(\alpha;\xi)$, see  Theorem 3.2 in \cite{BGR}. This change of determination is 
  a product of $k(n-k)-1$ {\it distinct} factors of degree 1,    
  to be compared with   (\ref{jumpGLS}-\ref{myster}).

 \section{Explicit value for low $n$}
 \label{lowvalues}
 \subsection{$n=2$}
 In that case, the PDF and the associated $\CI$ functions are easily determined. Let $\alpha=(\alpha_1, -\alpha_1)$, then 
 $$A=U(\theta) \diag (\alpha_1, -\alpha_1) U^\dagger(\theta),\qquad (0\le \theta \le \pi)\qquad
 \xi_1=A_{11} = \alpha_1 \cos \theta\,,$$ whence a support $\xi_1\in (-\alpha_1, \alpha_1)$ with a density 
 $$p(\xi | \alpha)= \oh \frac{\sin \theta \, d\theta}{d \xi_1}= \inv{2\alpha_1}\,, \qquad \CI=\oh\,.$$  
 As expected, the functions $p$ and $\CI$ are constant and discontinuous at the edges of their support. 
\\
 
 \noindent {\it Remark.} The parallel computation in the case of a real symmetric matrix with action of $\SO(2)$ leads 
 to a density singular on the edges of its support: 
 $p= \inv{\pi \sqrt{\alpha_1^2-\xi_1^2}}$.

  \subsection{$n=3$}

\begin{figure}[htb]
\begin{center}
\includegraphics[width=0.4\textwidth]{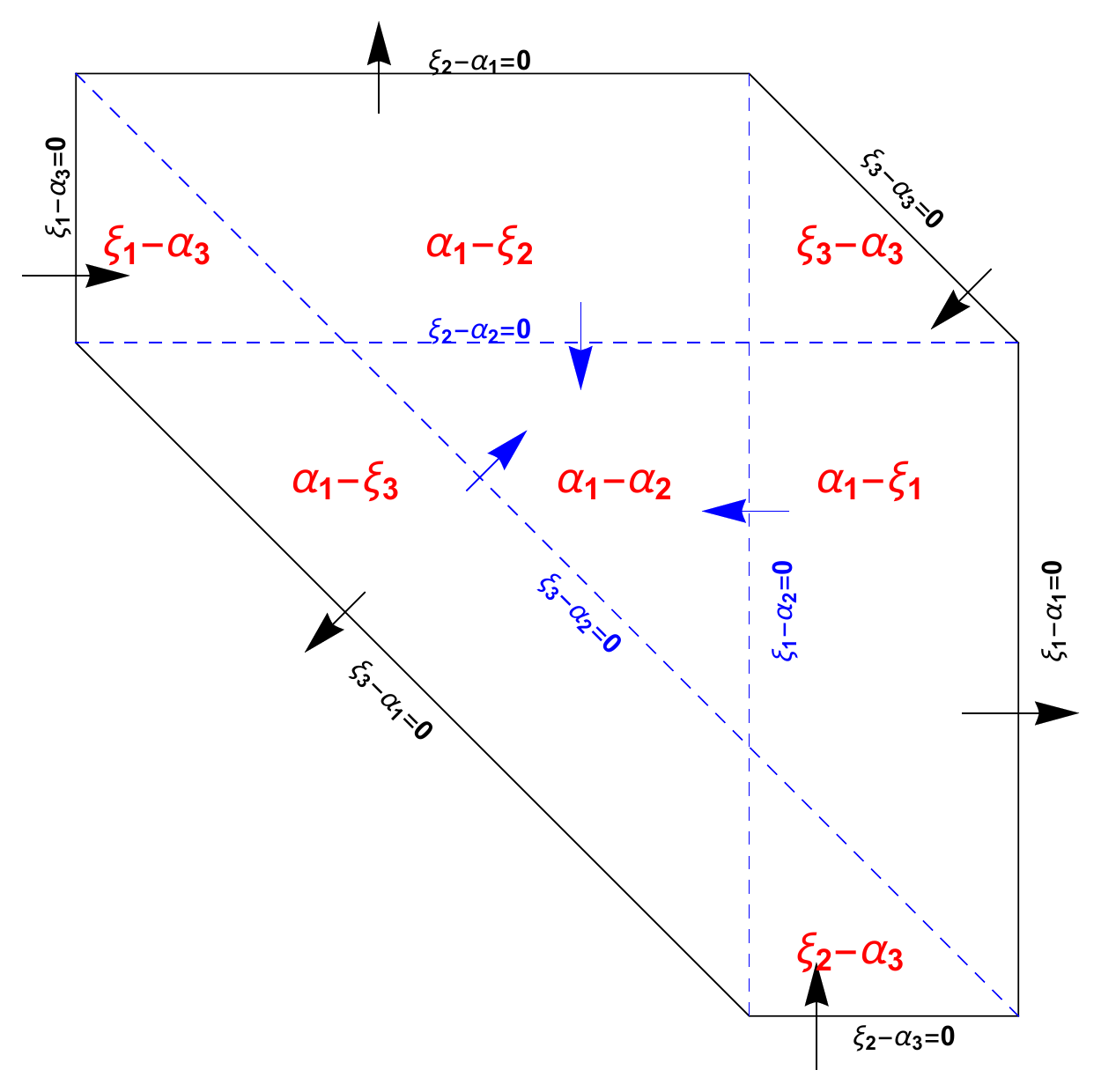}\qquad \includegraphics[width=0.4\textwidth]{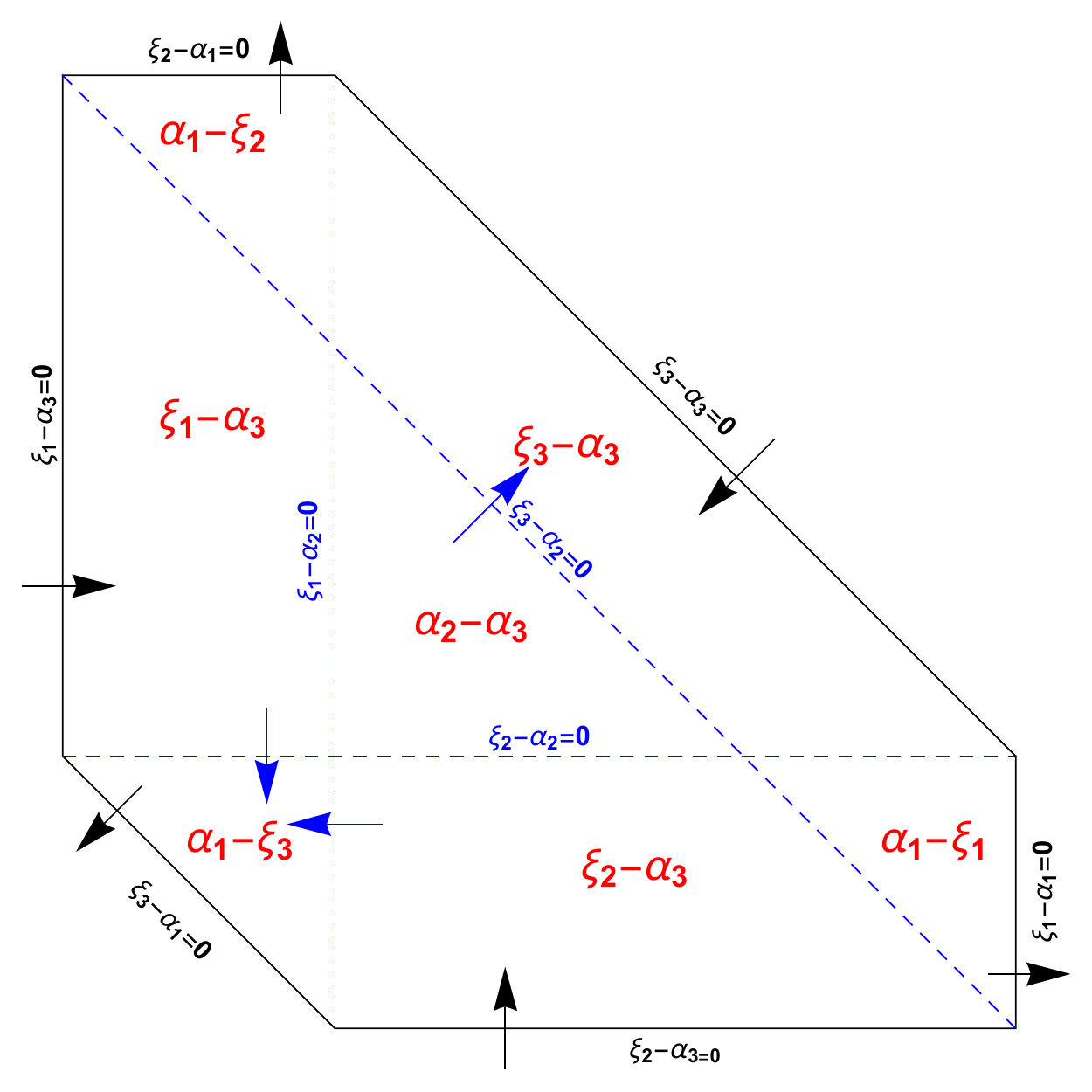}
\caption{Piecewise polynomial determination of the function $\CI_{su(3)}:= \mathrm{PDF} \times \Delta(\alpha)$
of the $\xi_i:=A_{ii}$, for SU(3). We assume $\alpha_1> \alpha_2>\alpha_	3$. 
Left: $\alpha_2>\oh(\alpha_1+\alpha_3)$. Right:   $\alpha_2<\oh(\alpha_1+\alpha_3)$.  }   
\label{Schur3-determinations}
\end{center}
\end{figure}

 For $n=3$, the function $\CI$  is readily computed  and given by the function displayed 
 in red on Fig.~\ref{Schur3-determinations}. It is normalized by $\int d^2 \xi\, \CI(\alpha;\xi)= \Delta(\alpha)/\Delta(\rho)$ as it should,
 see (\ref{p-CI}). Equivalent formulae have been given in \cite{GLS}.

There is an alternative way of presenting this result. Instead of giving the value of $\CI$ in each cell, start from the value 0 
at the exterior of the permutahedron and give the rule for the change of polynomial determination as one crosses a wall.
The rule is as follows: if the equation of a wall is $\xi_i -\alpha_j =0$, the function $\CI$ is incremented by $ \xi_i -\alpha_j$ when the wall 
is crossed in the direction of the arrow, see Fig.~\ref{Schur3-determinations}.
(Here one should remember that  $\sum_{i=1}^3 \xi_i = \sum_{i=1}^3 \alpha_i$.)

We conclude that in the case $n=3$, the changes of determination are by affine functions of $\xi$ that vanish on the 
hyperplanes (here lines) of singularity, in accordance with the $C^0$ class of differentiability. 

\medskip
 Remark. From the constraints on the honeycomb parameter, see below sect. \ref{polytope}, eq. (\ref{poly3}), 
 we may derive another expression of the same function 
 \bea\nonumber 
 \CI_{su(3)}(\alpha;\xi)&=& \min(\alpha_1, \alpha_1+\alpha_2-\xi_1)-\max(\alpha_2,\xi_2,\xi_3,\alpha_1+\alpha_3-\xi_1)\\
  \label{vol3}  &=& \min(\alpha_1-\alpha_2, \alpha_2-\alpha_3, \alpha_1-\xi_i, \xi_i-\alpha_3)\qquad i=1,2,3 \,, \eea
  which is manifestly non negative and fully symmetric in the $\xi$'s.

  Consider now two weights of SU(3), whose difference belongs to the root lattice.  As explained above, we denote them either by the 2-dimensional 
   vector of their (integral) Dynkin components $\lambda=\{\lambda_1,\lambda_2\}$ and $\delta=\{\delta_1,\delta_2\}$, or equivalently by their partition (Young) components $\alpha$ and  $\xi$. 
Then the number of integer points in the interval (\ref{intervalx22}) is the 
Kostka multiplicity, denoted   $\mult_\alpha(\xi)$  or $\mult_\lambda(\delta)$, and  reads
  \be\label{mult3} \mult_\alpha(\xi)   = \CI_{su(3)}(\alpha;\xi)+1  \,. \ee
  Expressed in Dynkin coordinates, this reads, if $(\delta_2 - \delta_1) \geq (\lambda_2 - \lambda_1)$ 
 \be
\mult_\lambda(\delta) =
 1 + \min[\frac{1}{3} ( 2 \lambda_1 + \lambda_2+\delta_1 - \delta_2 ), \lambda_2, 
  \frac{1}{3} (\lambda_1 + 2 \lambda_2+2 \delta_1 + \delta_2 ), \frac{1}{3} ( \lambda_1 + 2 \lambda_2-\delta_1 - 2 \delta_2)]
\ee
Otherwise, one should replace  $(\lambda, \delta)$ by the conjugate pair $(\overline \lambda, \overline \delta)$ in the above expression and use the fact that $\mult_{\overline\lambda}(\overline\delta)=\mult_\lambda(\delta)$.

These expressions are equivalent to those given in \cite{BGR}, sect. 7.1, for the multiplicities.

   \subsection{$n=4$}
  
\begin{figure}[htb]
\begin{center}
\includegraphics[width=0.7\textwidth]{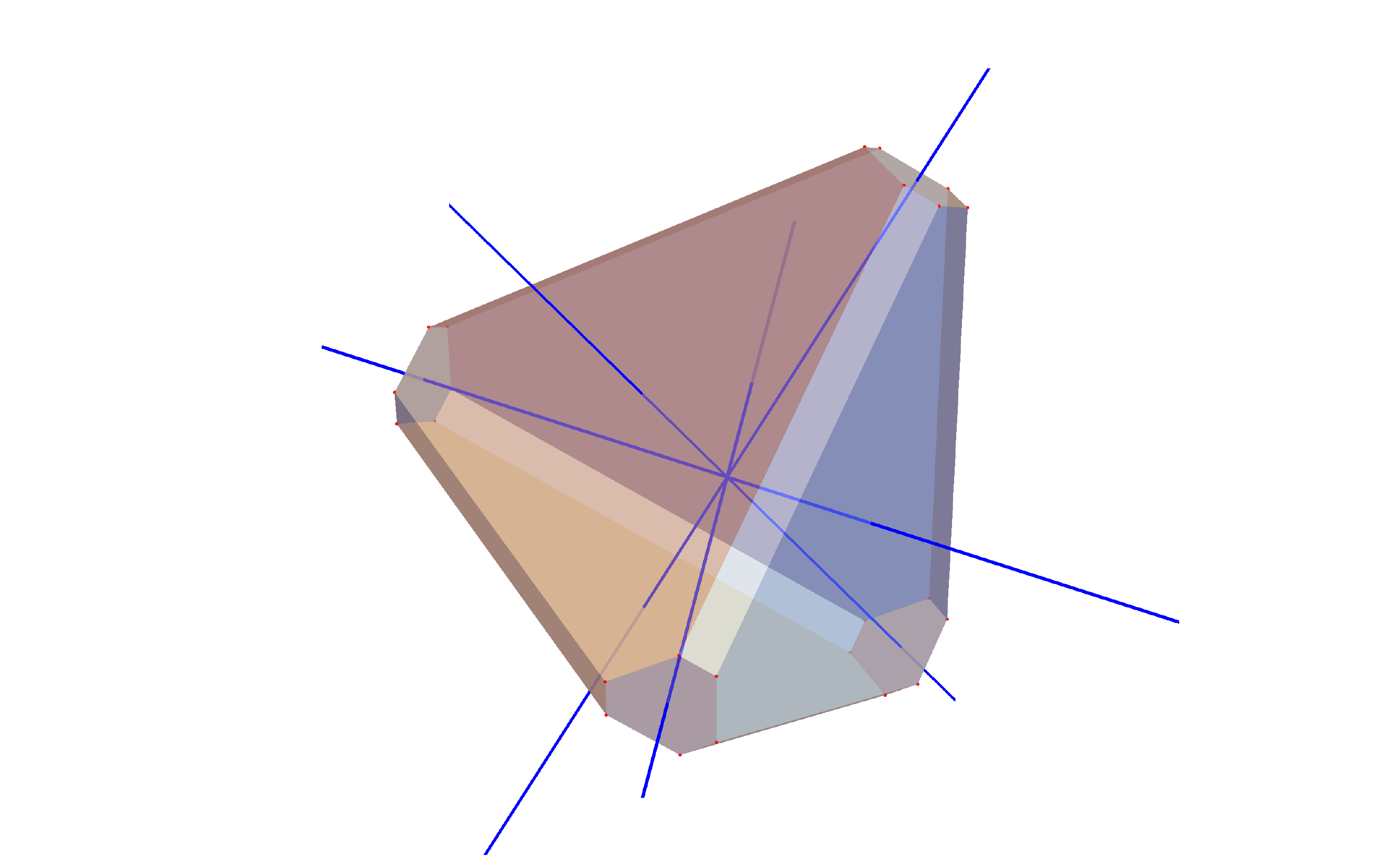}
\caption{The permutahedron for $n=4$ and $\alpha=(5,4,2,-11)$. \\The four coordinates
$\xi_i$ run between -11 and 5 along the four blue axes.}   
\label{jolipolytope}
\end{center}
\end{figure}

The permutahedron may be regarded as the convex polytope of points $\xi$ satisfying the inequalities
$$ \alpha_4\le \xi_i \le \alpha_1 \,,\  1\le i\le 4\,,  \qquad  \alpha_3+\alpha_4\le   \xi_i+\xi_j\le \alpha_1+\alpha_2,\ 1\le i<j\le 4\,.$$
It has thus four pairs of hexagonal faces $\xi_i= \alpha_1 \ \mathrm{or} \  \alpha_4$, $i=1,\cdots, 4$; and 
3 pairs of rectangular faces $\xi_i+\xi_j = \alpha_1+\alpha_2\  \mathrm{or}\   \alpha_3+\alpha_4$, see Fig.~\ref{jolipolytope}. (Note that 
$\xi_i+\xi_j =\alpha_1+\alpha_2 \Leftrightarrow \xi_k+\xi_\ell = \alpha_3+\alpha_4$, $i,j,k,\ell$ all distinct, 
since $\sum_i \xi_i=\sum_i \alpha_i=\tr A$.)
Moreover, we expect the function to be piecewise polynomial, with possible changes of determination across 
the hyperplanes of equation 
\\ \centerline{{\bf (i)}: $\xi_i= \alpha_2\ \mathrm{or} \  \alpha_3$;}\\ or \\  \centerline{{\bf (ii)}:
$\xi_i+\xi_j = \alpha_1+\alpha_3\  \mathrm{or}\   \alpha_1+\alpha_4$.} \\ These hyperplanes yield a partition 
of the permutahedron into (open) polyhedral cells.

From the expressions given in \cite{Z1,CZ1}  for the Horn volume, and taking the limit $\beta, \gamma \gg \alpha \approx \xi=\gamma-\beta$,
one finds (notation $A_i:= \alpha_{w(i)}-\xi_i, \ A_{ij}=A_i+A_j$, $A_{123}=A_1+A_2+A_3$)
\bea \CI_{su(4)}(\alpha;\xi)&=&\sum_{w\in \mathfrak{S}(4)} \epsilon(w) 
\varepsilon(A_1) 
\\ \nonumber &&
\bigg(\varepsilon(A_2) \Big(\frac{1}{6} \big(\left| A_{123}\right| ^3-\left|  A_{13}\right| ^3+\left| A_{23}\right| ^3-\left| A_3\right| ^3\big)   -\frac{1}{2} A_2 (A_{123} \left| A_{123}\right| + A_3 \left|A_3 \right| )\Big)
   \\ \nonumber &&
   +\varepsilon(A_{12}) \Big(\frac{1}{2} A_{12} (A_{123} \left|
 A_{123}\right| +A_3 \left|A_3 \right| )+\frac{1}{3} (\left|A_3\right| ^3-\left| A_{123}\right| ^3)\Big)\bigg)\,,  \qquad 
   \eea
 where $\epsilon(w)$ is the signature of $w$ and $\varepsilon (\cdot)$ is the sign function.\\

In principle, there is an alternative expression of $\CI_{su(4)}$, though not explicit,
 coming from its interpretation as the  volume of a 3d polytope, see sect. \ref{polytope}.  

\medskip
 Finally we have yet another expression, which makes explicit the   location of the singular (hyper)planes and the piecewise polynomial 
 determinations. This will be described now. 
  
In accordance with the results of \cite{GLS},  see above (\ref{jumpGLS}), we expect that the change of polynomial determination of $\CI_{su(4)}$ is (at least) quadratic across the hyperplanes of type (i),
and cubic across those of type (ii). This applies also to the vanishing of the function on the external faces,
with a quadratic, resp. cubic behavior on hexagonal, resp. rectangular faces. 
This is confirmed by the detailed and explicit expression of the jumps of $\CI_{su(4)}$ that we discuss now. 

\bigskip

In the vicinity of an hexagonal face, internal or external, $\xi_i=\alpha_j\,, i,j=1,\cdots,4$,  
$\CI_{su(4)}$ undergoes a change of determination $\Delta\CI_{su(4)}(\xi)$ of the form
\be\label{outerfaceC1}\Delta\CI_{su(4)}(\xi)= \oh  (\xi_i- \alpha_j)^2  \,  p_c(\xi)\ee
with $p_c(\xi)$ a degree 1 polynomial. 
 $\Delta\CI_{su(4)}(\xi)$ thus vanishes with its first order derivatives on that face, hence $\CI_{su(4)}$ is of class $C^1$. 
More precisely,  
as $\xi_1-\alpha_3$ increases through 0, for instance, which we denote by $\xi_1=\!\!\!\!\!\nearrow  \alpha_3$, 
$\CI$ is incremented by
\bea \label{delta13aa}
\nonumber &&\hskip69mm \ \quad \xi_2-\alpha_2 \qquad  \xi_3-\alpha_2 \qquad  \xi_4-\alpha_2 \qquad \\
\!\!\!\!\!\!\!\! \!\!\!\!\!\!\!\!  \Delta\CI_{su(4)}(\xi)\Big|_{\xi_1=\!\!\!\!\!\nearrow  \alpha_3}
\!\!\!\!\!\!\!\! &=&\!\!\!\!
 -\inv{6}(\xi_1-\alpha_3)^2 \times \begin{cases}  
3 (\xi_2-\alpha_4) +(\xi_1-\alpha_3) & \mathrm{if}\quad <0 \ \quad \qquad \ >0\ \quad \qquad\ >0\\
3 (\xi_3-\alpha_4)+(\xi_1-\alpha_3)  & \mathrm{if}\quad >0 \ \quad \qquad \ <0\ \quad \qquad\ >0\\
3 (\xi_4-\alpha_4)+(\xi_1-\alpha_3) & \mathrm{if}\quad >0 \ \quad \qquad \ >0\ \quad \qquad\ <0\\
3 (\alpha_2-\alpha_4) & \mathrm{if}\quad >0 \ \quad \qquad \ >0\ \quad \qquad\ >0\\
3 (\alpha_1-\alpha_2)  & \mathrm{if}\quad <0 \ \quad \qquad \ <0\ \quad \qquad\ <0\\
3 (\alpha_1-\xi_4)-(\xi_1-\alpha_3) & \mathrm{if}\quad <0 \ \quad \qquad \ <0\ \quad \qquad\ >0\\
3 (\alpha_1-\xi_3)-(\xi_1-\alpha_3) & \mathrm{if}\quad <0 \ \quad \qquad \ >0\ \quad \qquad\ <0\\
3 (\alpha_1-\xi_2)-(\xi_1-\alpha_3) & \mathrm{if}\quad >0 \ \quad \qquad \ <0\ \quad \qquad\ <0
\end{cases}
\eea
or, in a more compact form,
\be \label{delta13ab} \Delta\CI_{su(4)}(\xi)\Big|_{\xi_1=\!\!\!\!\!\nearrow  \alpha_3}=-\inv{2}(\xi_1-\alpha_3)^2 \Big(\CI_{su(3)}((\alpha_1,\alpha_2,\alpha_4);(\xi_2,\xi_3,\xi_4))
+\frac{\eta_{13}}{3} (\xi_1-\alpha_3)\Big) \,,
 \ee
   in terms of the function defined above in Fig. \ref{Schur3-determinations} or in eq. (\ref{vol3}), and 
where $\eta_{13}=+1,0,-1$ depending on the case, as read off the expression (\ref{delta13aa}): $\eta_{13}=1$ on the first three lines of (\ref{delta13aa}),
then $0$ on the next two, and $-1$ on the last three.
In other words, if  
\be  
h((\beta_1,\beta_2,\beta_3); (\zeta_1,\zeta_2,\zeta_3)):= 
\begin{cases} +1 &\mathrm{if}\ \sign(\zeta_1-\beta_2)+\sign(\zeta_{ 2}-\beta_2)+\sign(\zeta_{ 3}-\beta_2) =+1\\\
0 &\mathrm{if}\ \sign(\zeta_1-\beta_2)+\sign(\zeta_{ 2}-\beta_2)+\sign(\zeta_{ 3}-\beta_2) =\pm 3 \\
-1 &\mathrm{if}\ \sign(\zeta_1-\beta_2)+\sign(\zeta_{ 2}-\beta_2)+\sign(\zeta_{ 3}-\beta_2)  =-1
\end{cases} \,,\ee
then  $\eta_{13}= h((\alpha_1,\alpha_2,\alpha_4),(\xi_2,\xi_3,\xi_4))$, or more generally, 
$\eta_{ij}= h((\alpha_1,\cdots,\widehat{\alpha_j},\cdots,\alpha_4),(\xi_1,\cdots,\widehat{\xi_i},\cdots,\xi_4))$.
(As usual, the caret means omission.)

The analogous increments through $\xi_1=\alpha_2$ are given by 
similar expressions, where $\alpha_2$ and $\alpha_3$ have been swapped, 
{\bf and} the overall sign changed, thus 
\be \label{delta12ab} \Delta\CI_{su(4)}(\xi)\Big|_{\xi_1=\!\!\!\!\!\nearrow  \alpha_2}=\inv{2}(\xi_1-\alpha_2)^2 \Big(\CI_{su(3)}((\alpha_1,\alpha_3,\alpha_4);(\xi_2,\xi_3,\xi_4))
+\frac{\eta_{12}}{3} (\xi_1-\alpha_2)\Big) \,,
 \ee
with now $\eta_{12}$ determined by the signs of the three differences $\xi_i-\alpha_3$, $i=2,3,4$.

These formulae also apply to the external hexagonal faces, for example $\xi_1=\alpha_4$, but now $\CI$ vanishes 
on one side, thus the formula actually gives the value of $\CI$:
\be\label{exter-face}  \Delta\CI_{su(4)}(\xi)\Big|_{\xi_1=\!\!\!\!\!\nearrow  \alpha_4}=\CI_{su(4)}(\xi)\Big|_{\xi_1= \alpha_4+0}=
\inv{2}(\xi_1-\alpha_4)^2 \Big(\CI_{su(3)}((\alpha_1,\alpha_2,\alpha_3);(\xi_2,\xi_3,\xi_4))
+\frac{\eta_{14}}{3} (\xi_1-\alpha_4)\Big) \,,
 \ee
where $\eta_{14}$ determined by the signs of the three differences $\xi_i-\alpha_2$, $i=2,3,4$.
Note that the overall sign is in agreement with the positivity of the functions $\CI_{su(4)}$  and $\CI_{su(3)}$.
Likewise, across the ``upper" face $\xi_1=\alpha_1$,
\be  \Delta\CI_{su(4)}(\xi)\Big|_{\xi_1=\!\!\!\!\!\nearrow  \alpha_1}=-\CI_{su(4)}(\xi)\Big|_{\xi_1= \alpha_1-0}=
-\inv{2}(\xi_1-\alpha_1)^2 \Big(\CI_{su(3)}((\alpha_2,\alpha_3,\alpha_4);(\xi_2,\xi_3,\xi_4))
+\frac{\eta_{11}}{3} (\xi_1-\alpha_1)\Big) \,.
 \ee
Finally, since  $\CI_{su(4)}$ is a symmetric function of the $\xi_i$, $i=1,\cdots,4$, 
the changes of determination across walls of equation $\xi_{2,3, 4}=\alpha_{j}$  are obtained from the previous expressions 
by a  permutation of the $\xi$'s. 
\bigskip

We now turn to the rectangular faces. Across a rectangular face of equation  $\xi_i+\xi_j=\alpha_k+\alpha_l$, $\CI_{su(4)}$ is of class 
$C^2$, hence its change of determination $\Delta\CI_{su(4)}$ vanishes cubically. 
We find 
$$\Delta\CI_{su(4)}\Big|_{\xi_i+\xi_j=\!\!\!\!\!\nearrow  \alpha_k+\alpha_l}=
\pm \inv{6} \big(\xi_i+\xi_j-\alpha_k-\alpha_l\big)^3$$
The overall sign is determined by the positivity for external faces  $\xi_i+\xi_j=\alpha_1+\alpha_2$ or  $\xi_i+\xi_j=\alpha_3+\alpha_4$
\bea\label{outerfaceC2}
\label{exter-face12}\Delta\CI_{su(4)}\Big|_{\xi_i+\xi_j=\!\!\!\!\!\nearrow  \alpha_1+\alpha_2}=&
-\CI(\xi)\Big|_{\xi_i+\xi_j=\alpha_1+\alpha_2-0}&= \frac{1}{6}(\xi_i+\xi_j-\alpha_1-\alpha_2)^3\,,\\
\label{exter-face34} \Delta\CI_{su(4)}\Big|_{\xi_i+\xi_j=\!\!\!\!\!\nearrow  \alpha_3+\alpha_4}=&
\CI(\xi)\Big|_{\xi_i+\xi_j=\alpha_3+\alpha_4+0}&=    \frac{1}{6}(\xi_i+\xi_j-\alpha_3-\alpha_4)^3\,. \eea
Internal walls of type (ii) may be written as $\xi_i+\xi_j =\alpha_1+\alpha_3$ or $\xi_i+\xi_j =\alpha_2+\alpha_3$, with $1\le i< j\le 4$.
As $\xi_i+\xi_j-\alpha_{1\atop 2}-\alpha_3$ increases through 0, 
the {\it change} of polynomial determination is  given by 
\bea 
\label{casea}\Delta\CI_{su(4)}\Big|_{\xi_i+\xi_j=\!\!\!\!\!\nearrow  \alpha_1+\alpha_3}&=& -\inv{6}(\xi_i +\xi_j-\alpha_1-\alpha_3)^3 \\
\label{caseb}
\Delta\CI_{su(4)}\Big|_{\xi_i+\xi_j=\!\!\!\!\!\nearrow  \alpha_2+\alpha_3}&=&
 \begin{cases}
\inv{6}(\xi_i +\xi_j-\alpha_2-\alpha_3)^3  & \mathrm{if}\  \alpha_3 \le  \xi_i,\xi_j\le \alpha_2  \\
0 & \mathrm{otherwise}\,.\end{cases}\eea

\begin{figure}[htb]
\begin{center}
\includegraphics[width=0.45\textwidth]{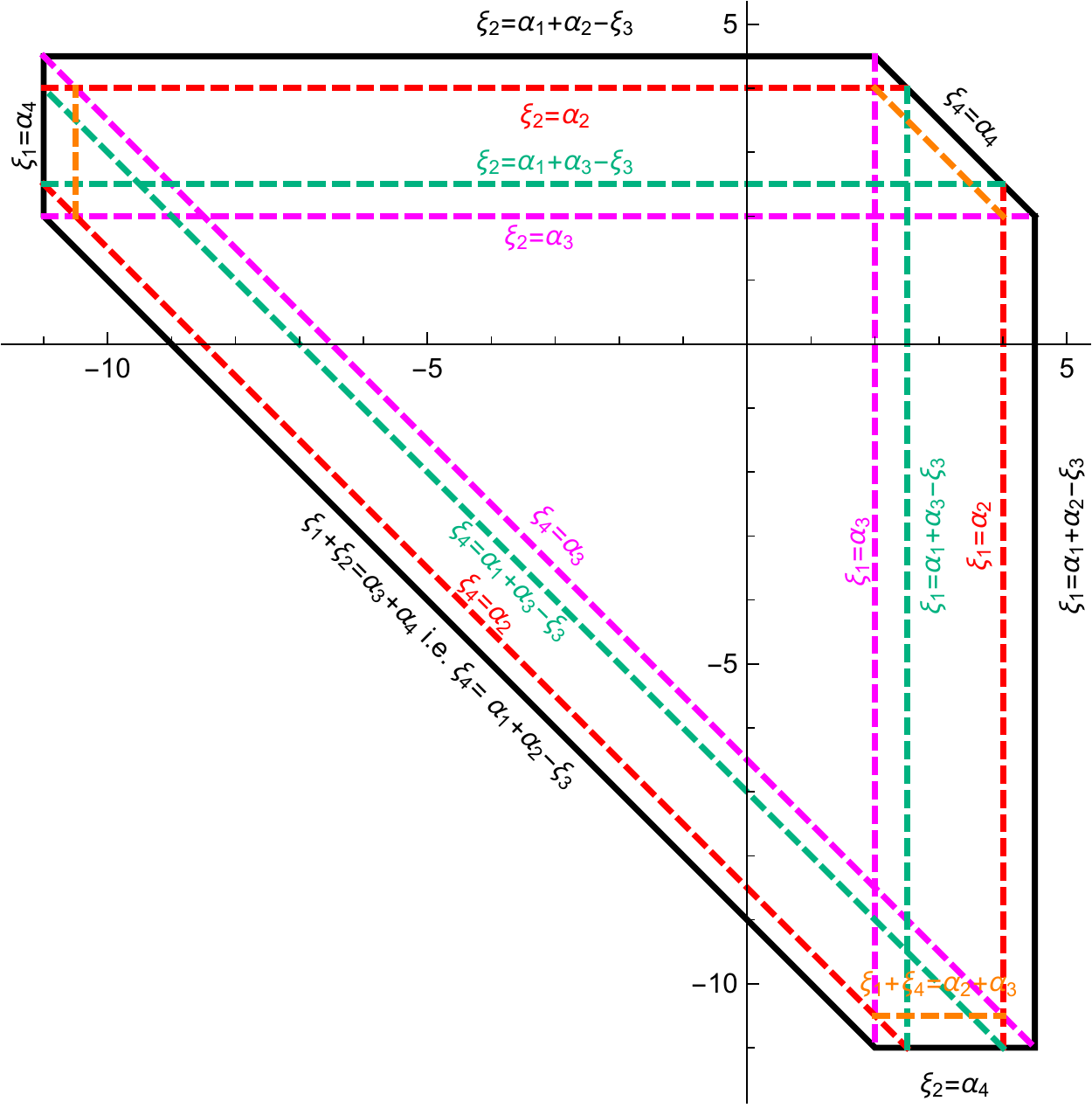}
\caption{ SU(4): A cross-section of the permutahedron, again for for $\alpha=\{5,4,2,-11\}$, along the plane at fixed 
$\xi_3 = (\alpha_1+\alpha_2)/2$.  (Recall that $\xi_4=\sum_i\alpha_i -\xi_1-\xi_2-\xi_3$.) 
There are $ 1+3\times (9+3) =37$ cells of polynomiality. 
Note that the orange lines 
(of equation $\xi_i+\xi_j=\alpha_2+\alpha_3$, $i,j=1,2,4$) do not extend to the boundary of the permutahedron, reflecting the condition 
(\ref{caseb}).  
}
\label{cross-sections-n=4}  
\end{center}
\end{figure}

Example. 
Fig. \ref{cross-sections-n=4}  displays a cross-section of the permutahedron for $\alpha=\{5,4,2,-11\}$ at 
$\xi_3=(\alpha_1+\alpha_2)/2$.
The hyperplanes $\xi_i+\xi_j=\alpha_1+\alpha_3$, resp. $=\alpha_2+\alpha_3$ intersect these cross-sections along the green 
lines, resp. the orange lines. 
There are $ 1+3\times (9+3) =37$ cells of polynomiality for $\max(\alpha_2, (\alpha_1-\alpha_2+\alpha_3) )\le \xi_3\le \alpha_1$.

\medskip
Remarks.\\
1. The above expressions have been obtained in a semi-empirical way, checking on many cases
their consistency with 
the original expression (\ref{CI}). A direct and systematic proof would clearly be desirable. \\
2. Denote by $\{i,j,k,\ell\}$ a permutation of $\{1,2,3,4\}$.  
The  first case (\ref{casea}) above,   ({green}   lines on Fig. \ref{cross-sections-n=4}),
occurs only if $\alpha_3\le \xi_i, \xi_j\le \alpha_1$ and $\alpha_4\le \xi_k,\xi_\ell \le \alpha_2$.
Indeed if $\xi_i+\xi_j=\alpha_1+\alpha_3 \Leftrightarrow \xi_k+\xi_\ell =\alpha_2+\alpha_4$, then $\xi_i=\underbrace{(\alpha_1-\xi_j)}_{\ge 0} +\alpha_3 \ge \alpha_3$ and 
of course $\xi_i\le \alpha_1$, while $\xi_k= \alpha_2+ \underbrace{(\alpha_4-\xi_\ell)}_{\le 0} \le \alpha_2$ and of course $\xi_k\ge \alpha_4$. In other 
words, the planes  $\xi_i+\xi_j=\alpha_1+\alpha_3$ do not intersect the permutahedron if those conditions are not fulfilled.\\
3. In contrast, in the second case  (\ref{caseb})  ({orange}  
lines on Fig. \ref{cross-sections-n=4}), there is
a change of polynomial determination only if $\alpha_3\le \xi_i,\xi_j\le \alpha_2$ although the planes $\xi_i+\xi_j=\alpha_2+\alpha_3$ 
intersect the permutahedron irrespective of that condition.\\
This illustrates the well known fact that the forms given above in (\ref{equsing}) for the loci of change of polynomial determinations are
only necessary conditions. It may be that the function $\CI$ is actually regular across some of these (hyper)planes.

\section{Relation of $\CI$ with multiplicities}   
\label{sec4}

\subsection{The Kostant multiplicity formula }   
\label{kostantmult}

The multiplicity of the weight $\delta$ in the 
 irreducible representation of highest weight $\lambda$, aka the Kostka number, 
  may be written in various ways, e.g. following Kostant  \cite{Kostant:multiplicity}
\be\label{kost-partition-formula}\mult_\lambda(\delta) =\sum_{w\in W} \epsilon_w P(w(\lambda+\rho) -\delta-\rho)\ee
 where $P$ is Kostant's partition function (\cite{Hall}, Theorem 10.29). 
 $P(\Gb)$ gives the number of ways an element $\Gb$ of the root lattice may be decomposed as a non-negative integer linear combination  of positive roots. 
 
 On the other hand, it has been pointed out by Heckman (\cite{He82}, see also \cite{GLS}) that {\it asymptotically}, for large weights $\lambda$ and $\delta$,
one may write a semi-classical approximation of $\mult_\lambda(\delta)$. We shall recover below this result from the more general formula
(\ref{asy-scaling}). 

\bigskip

The Kostant partition functions for rank 2 Lie algebras are given in \cite{Tarski:Kostant} and in \cite{Capparelli:Kostant}.
The first reference also gives the partition function for $A_3$. 
The expression obtained by \cite{Capparelli:Kostant} for $B_2$ is more compact and this is the one that we give below.
We used these formulae to check  the consistency of expressions obtained by other means, 
 using  honeycombs or other pictographs, or counting the integer points in Berenstein-Zelevinsky (BZ) polytopes, see below.


If $\kappa$ is a weight, with components $k_i$ in the basis of simple roots (``Kac labels'') (warning: {\sl not} the basis of fundamental weights (``Dynkin labels'')!), $P(\kappa) $ is given by: \\
{{\small $\SU(3)$: \; 
$P(k_1,k_2) =\begin{cases}  \min ({k_1},{k_2})+1 & \text{if}\  {k_1}\ge 0 \;  \mathrm{and}\; {k_2}\ge 0 \cr
0& \text{otherwise} 
 \end{cases}$
 
\smallskip
\noindent
{\small $\SU(4)$: \;  $P(k_1,k_2,k_3) =$ \\
$
\left(
\begin{array}{cc}
 0 & {k_1}<0\lor {k_2}<0\lor {k_3}<0 \\
 \frac{1}{6} ({k_2}+1) ({k_2}+2) ({k_2}+3) & {k_2}\leq {k_1}\land {k_2}\leq {k_3} \\
 \frac{1}{6} ({k_1}+1) ({k_1}+2) (-2 {k_1}+3 {k_2}+3) & {k_1}\leq {k_2}\leq {k_3} \\
 \frac{1}{6} ({k_1}+1) ({k_1}+2) (-{k_1}+3 {k_3}+3) & {k_1}\leq {k_3}\leq {k_1}+{k_3}\leq {k_2} \\
 \frac{1}{6} ({k_1}+1) ({k_1}+2) (-{k_1}+3 {k_3}+3)- { {k_1 - k_2 +k_3+ 2}\choose{ 3 } }
    & {k_1}\leq {k_3}\leq {k_2}\leq {k_1}+{k_3} \\
 \frac{1}{6} ({k_3}+1) ({k_3}+2) (3 {k_2}-2 {k_3}+3) & {k_3}\leq {k_2}\leq {k_1} \\
 \frac{1}{6} ({k_3}+1) ({k_3}+2) (3 {k_1}-{k_3}+3) & {k_3}\leq {k_1}\leq {k_1}+{k_3}\leq {k_2} \\
 \frac{1}{6} ({k_3}+1) ({k_3}+2) (3 {k_1}-{k_3}+3)- { {k_1 - k_2 +k_3+ 2}\choose{ 3 } }
   & {k_3}\leq {k_1}\leq {k_2}\leq {k_1}+{k_3} \\
\end{array}
\right).
$}
\\
 Notice that for $\SU(4)$ there are seven cases (Kostant chambers of polynomiality) and that the last three cases are obtained from the previous three by exchanging $k_1$ and $k_3$.

\smallskip
\noindent
{\small $B_2$:\;  $P(k_1,k_2) =$
$
\left(
\begin{array}{cc}
 0 & {k_1}<0\lor {k_2}<0 \\
 {b}({k_2}) & {k_2}\leq {k_1} \\
 \frac{1}{2} ({k_1}+1) ({k_1}+2) & 2 {k_1}\leq {k_2} \\
 \frac{1}{2} ({k_1}+1) ({k_1}+2)-{b}(2 {k_1}-{k_2}-1) & {k_1}\leq {k_2}\leq 2 {k_1} \\
\end{array}
\right)
$}
\\
Here the function $b$ is defined by  $b(x)= \left(-\lfloor \frac{x+1}{2}\rfloor+x+1\right) \left(\lfloor\frac{x+1}{2}\rfloor+1\right)$, in terms of 
the integer part function.}


 \subsection{The $\CI$--multiplicity relation }   
 \label{CI-mult}
We follow here the same steps as in the case of the relation betweeen the  Horn volume and LR coefficients~\cite{CMSZ1, CMSZ2}.
Take  $\lambda'=\lambda+\rho$, where $\lambda$ is the h.w. of the irrep $V_\lambda$ and $\rho$ is the Weyl vector.
\bea\label{I-mult-rel}
\CI(\lambda';\delta)&=& \frac{\Delta(\lambda')}{\Delta(\rho)} \int_{\R^{n-1}}\frac{dx}{(2\pi)^{n-1}}  \CH(\lambda';\delta)  e^{-\ii \langle x,\delta\rangle}\\
\nonumber &=&  \int_{\R^{n-1}}\frac{dx}{(2\pi)^{n-1}} \underbrace{ \CH(\lambda';\delta)  \dim V_\lambda}_{\frac{\hat\Delta(e^{\ii x})}{\Delta(\ii x)} \chi_\lambda(e^{\ii x})}  e^{-\ii \langle x,\delta\rangle}\\
\nonumber &=&   \int_{\R^{n-1}/2\pi P^\vee}\frac{dx}{(2\pi)^{n-1}} \sum_{\psi\in 2\pi P^\vee}
\frac{\hat\Delta(e^{\ii (x+\psi)})}{\Delta(\ii (x+\psi))} \chi_\lambda(e^{\ii (x+\psi)})  e^{-\ii \langle (x+\psi),\delta\rangle}
\\
\nonumber &=&    \int_{\R^{n-1}/2\pi P^\vee}\frac{dx}{(2\pi)^{n-1}} \sum_{\psi\in 2\pi P^\vee}  e^{\ii \langle  \lambda+\rho-
\delta, \psi\rangle}
\frac{\hat\Delta(e^{\ii x})}{\Delta(\ii (x+\psi))} \chi_\lambda(e^{\ii x})  e^{-\ii \langle x,\delta\rangle}
 \eea
Following  Etingof and Rains \cite{ER} 
one writes 
\bea\label{ER}
\sum_{\psi\in 2\pi P^\vee} {e^{\ii \langle \rho, \psi\rangle}}
\frac{\hat\Delta(e^{\ii x})}{\Delta(\ii (x+\psi))}&=&\sum_{\kappa\in K} r_\kappa\chi_\kappa(e^{\ii x}) \\
\nonumber \sum_{\psi\in 2\pi P^\vee} 
\frac{\hat\Delta(e^{\ii x})}{\Delta(\ii (x+\psi))}&=&\sum_{\kappa\in \hat K} \hat r_\kappa\chi_\kappa(e^{\ii x}) \,,
\eea
where the finite sets of weights $K$ and $\hat K$ and the rational coefficients $r_\kappa,\,  \hat r_\kappa$ have been 
defined in \cite{CMSZ1, CMSZ2}, and one finds, with $Q$ the root lattice, 
\be\label{I-mult}
\CI(\lambda';\delta) =\begin{cases} \sum_{\kappa\in K} r_\kappa\, \sum_\tau C_{\lambda \kappa}^\tau  \, \mult_\tau(\delta) &\mathrm{if} \ 
\lambda-\delta\in Q \cr
\sum_{\kappa\in \hat K} \hat r_\kappa  \, \sum_\tau C_{\lambda \kappa}^\tau  \, \mult_\tau(\delta) &\mathrm{if}\  \lambda+\rho-\delta\in Q\,. 
\end{cases}
\ee
\\
For su(3), the two formulae boil down to the same simple expression, since $\rho\in Q$ and $K=\hat K=\{0\}$, hence
\be\label{CI-mult-su3} \CI_{su(3)}(\lambda';\delta)= \mult_\lambda(\delta)\,.\ee
This expression is compatible with that given in (\ref{mult3}), thanks to a peculiar identity that holds in su(3):
$\CI_{su(3)}(\lambda+\rho; \delta)= \CI_{su(3)}(\lambda;\delta)+1$  or equivalently $\mult_{\lambda+\rho}(\delta) = \mult_\lambda(\delta) +1\,.$
The latter is itself  obtained in the large $\mu,\nu=\mu+\delta$ limit of the more general
identity $C_{\lambda+\rho\, \mu+\rho}^{\nu+\rho}=C_{\lambda\mu}^\nu+1$ already mentioned in \cite{CZ0,CZ1}.
\\[8pt]
For su(4), we have two distinct relations
\be\label{CI-mult-su4} \CI_{su(4)}(\lambda';\delta) =\begin{cases}\inv{24} \( 9\, \mult_\lambda(\delta) +\sum_\tau  C_{\lambda\,\{1,0,1\}}^\tau\mult_{\tau}(\delta)  \)& \mathrm{if}\ \lambda-\delta \in Q \cr
\inv{6} \sum_\tau C_{\lambda\,\{0,1,0\}}^\tau\mult_{\tau}(\delta)& \mathrm{if}\ \lambda-\delta-\rho \in Q \,.
\end{cases}
\ee
In particular for $\lambda=0$, 
\be \CI_{su(4)}(\rho;\delta) =\begin{cases}\inv{24} \( 9\, {\delta}_{\delta 0} +\mult_{\{1,0,1\}}(\delta)  \)& \mathrm{if}\ \delta \in Q \cr
\inv{6}\mult_{\{0,1,0\}}(\delta)  & \mathrm{if}\ \delta-\rho \in Q \
\end{cases}
\ee
thus for $\delta=0$, $ \CI_{su(4)}(\rho;0) =\oh$ and for $\delta=\{1,0,1\}$ (or any of its Weyl images), $ \CI_{su(4)}(\rho;\delta) =\inv{24}$,
while for $\delta=\{0,1,0\}$ (or any of its Weyl images), $ \CI_{su(4)}(\rho;\delta) =\inv{6}$.

{\bf Remark}. Inverting the $\CI$--multiplicity formula  (\ref{I-mult})  is an interesting question that has been addressed in  \cite{CMS}, sect. 6.


\subsection{Polytopes and reduced \KT honeycombs}    
\label{polytope}
Recall that  Knutson--Tao (KT) honeycombs or other pictographs relevant for the LR coefficients of $su(n)$ depend on $(n-1)(n-2)/2$ parameters.
For example, in the \KT honeycombs,  the $3n$ external edges carry the components of $\alpha,\beta,\gamma$, while
each internal line carries a number, such that at each vertex the sum of the incident numbers vanishes. Moreover,
for each of the $3 n (n+1)/2$ internal edges, one writes a certain inequality between those numbers. 
 In the limit $\beta\sim \gamma \gg \alpha, \xi$, one third  of these inequalities is automatically satisfied and
one is left with $n(n+1)$ linear inequalities on the $(n-1)(n-2)/2$  parameters.

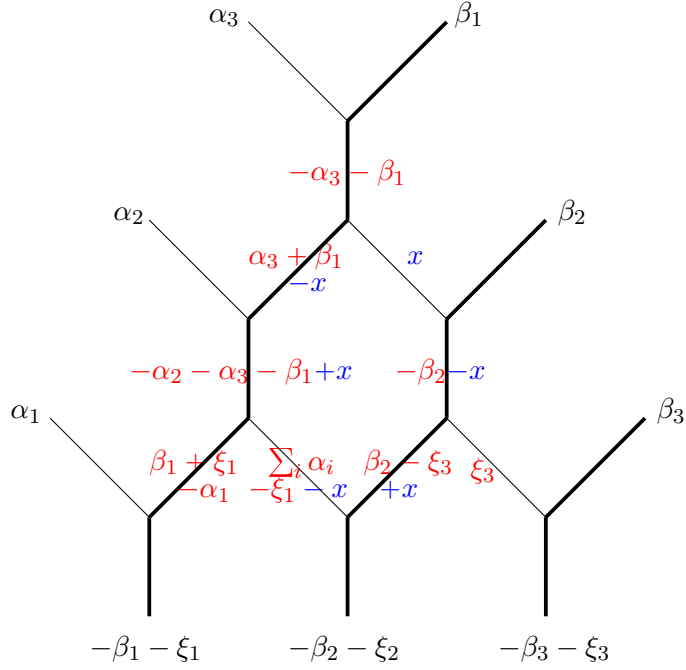
\begin{figure}[tb]
\setlength{\unitlength}{1.5pt}
{\begin{picture}(-400,60)
                           \put (115,45){$\alpha_3$}           { \linethickness{0.5mm}  \put(150,20){\line(1,1){25}}}
                            \put(150,20){\line(-1,1){25}} \put (177,45){$\beta_1$}
                    { \linethickness{0.5mm}      \put(150,20){\line(0,-1){25}} }   \put(135,5){$\red{ -\alpha_3-\beta_1}$}
                                        \put(125,-16){$\red{  \alpha_3+\beta_1 } $}
                                       \put(135,-23){$  \Blue{ - x} $}   
                                         \put(165,-16){${ \Blue{x}} $}                     
  \put (91,-5){$\alpha_2$}  \put(125,-30){\line(-1,1){25}} {\linethickness{0.5mm}\put(150,-5){\line(-1,-1){25}}}
\put(150,-5){\line(1,-1){25}}  { \linethickness{0.5mm}\put(175,-30){\line(1,1){25}}}
 			\put(95,-45){$\red{  -\alpha_2-\alpha_3-\beta_1} { \Blue{+ x}}$}\put(162,-45){$\red{  -\beta_2} {\Blue{-x}}$} 
                                                     { \linethickness{0.5mm} \put(125,-30){\line(0,-1){25}} \put(175,-30){\line(0,-1){25}} \put (203,-5){$\beta_2$}
 \put (65,-55){$\alpha_1$} } \put(75,-55){\line(1,-1){25}}   \put(100,-68){$\red{\beta_1+\xi_1}$}  \put(107,-75){$\red{-\alpha_1}$}  
  \put(130,-68){$\red{\sum_i \alpha_i}$}\put(125,-75){$ \red{-\xi_1}\Blue{ -x}$}
  \put(154,-68){$\red{\beta_2-\xi_3}$}  \put(158,-75){$\Blue{ +x}$} \put(181,-70){$\red{  \xi_3} $} 
 { \linethickness{0.5mm} \put(125,-55){\line(-1,-1){25}} } \put(125,-55){\line(1,-1){25}}{\linethickness{0.5mm}\put(175,-55){\line(-1,-1){25}}}
 \put(175,-55){\line(1,-1){25}}   \linethickness{0.5mm}{\put(225,-55){\line(-1,-1){25}}} \put (228,-55){$\beta_3$}
              \linethickness{0.5mm}   \put(100,-80){\line(0,-1){25}} \put(150,-80){\line(0,-1){25}}\put(200,-80){\line(0,-1){25}} 
                  \put(85,-115){$ -\beta_1-\xi_1$}         \put(135,-115){$ -\beta_2-\xi_2$}  \put(188,-115){$ -\beta_3-\xi_3$}       
 \end{picture}    \vskip6cm            
}
\caption{Knutson--Tao's honeycomb for $n=3$ for $\gamma,\beta\gg \alpha, \xi=\gamma-\beta$. Thick lines carry the large values $\beta_1>\beta_2>\beta_3\gg \alpha_i,\xi_i$.}
\label{HoneyComb-Schur-SU3}
\end{figure}

This is illustrated on
Fig. \ref{HoneyComb-Schur-SU3} for $n=3$. There, the choice of the parameters is such that the large numbers  
 of order O$(\beta,\gamma)$ are carried by North--South 
and NE--SO lines (heavy lines on the figure), while the NO--SE lines  carry numbers of order O$(\alpha,\xi)$.
Inequalities attached to the latter are automatically satisfied. 
\\[-30pt]
In general, the surviving inequalities are of the type $c\le a \le b $  for all patterns of the type
\hskip-14mm
{\begin{picture}(-300,55) \put (30,25){$a$}  \put(50,15){\line(-1,1){15}} \put(50,15){\linethickness{0.5mm} \line(1,1){15}}  \put(80,15){\line(-1,1){15}} \put (80,15){$c$}    \put(50,15){\linethickness{0.5mm}\line(0,-1){15}}  \put(50,0){\line(1,-1){15}} \put (65,-15){$b$} 
  \end{picture} 
 \\ within the honeycomb. \ \vskip1cm 

One then sees that the \KT honeycomb boils down to a Gelfand--Tsetlin (GT) triangle, 
see Fig. \ref{HoneyComb-GTtriangle}, and one is left with the  \GT inequalities
\bea \label{GTineq} \alpha_n \le x_{n-1}^{(n-1)} \le \alpha_{n-1}\le  &\cdots& \le x_1^{(n-1)}\le \alpha_1 \\
\nonumber      x_{i+1}^{(j+1)} &\le x_{i}^{(j)} \le& x_{i}^{(j+1)}\ ,\qquad {1 \le i,j \le n-2 }
\eea
together with the $n-2$ conservation laws
\be\label{conslaws} \sum_{i=1}^n \alpha_i = \xi_1 +\sum_{j=1}^{n-1} x_j^{(n-1)}= \xi_1+\xi_2+\sum_{j=1}^{n-2} x_j^{(n-2)}=\cdots =\sum_{i=1}^{n-2} \xi_i+x_1^{(2)}+x_2^{(2)}  = \sum_{i=1}^n \xi_i\,.\ee

{\small According to the well-known rules, the  semi-standard tableaux corresponding to that triangle must have 
$x_1^{(1)}=\xi_n$ boxes containing 1 (necessarily in the first row); $x_1^{(2)} + x_2^{(2)} -x_1^{(1)} =\xi_{n-1}$ boxes containing 2
(in the first two rows); etc;  and $\xi_1$ boxes containing $n$.}

\medskip
Relations (\ref{GTineq}-\ref{conslaws}) define a polytope $\CP(\alpha;\xi)$ in $\R^{(n-1)(n-2)/2}$, whose volume is the function $\CI_{su(n)}(\alpha;\xi)$,
as we show below in sect. \ref{Asymptotics}.

\begin{figure}[tb]
{\begin{picture}(-400,60)
                           \put (72,45){$\alpha_n$}           { \linethickness{0.5mm}  \put(110,20){\line(1,1){25}}}
                            \put(110,20){\line(-1,1){25}} \put (137,45){$\beta_1$}
                    { \linethickness{0.5mm}      \put(110,20){\line(0,-1){25}} } 
                    \put(110,-15){$\red{ x^{(n-1)}_{n-1}}$}  
                      \put (37,-5){$\alpha_{n-1}$}  \put(85,-30){\line(-1,1){25}} {\linethickness{0.5mm}\put(110,-5){\line(-1,-1){25}}}
\put(110,-5){\line(1,-1){25}}  { \linethickness{0.5mm}\put(135,-30){\line(1,1){25}}}
                                                     { \linethickness{0.5mm} \put(85,-30){\line(0,-1){25}} \put(135,-30){\line(0,-1){25}}
  \put (163,-5){$\beta_2$}
  
 \put (12,-55){$\alpha_{n-2}$} } \put(35,-55){\line(1,-1){25}} 
 \put(136,-70){$\red{ x^{(n-2)}_{n-2}}$} 
  \put(85,-70){$\red{ x^{(n-1)}_{n-2}}$}
 { \linethickness{0.5mm} \put(85,-55){\line(-1,-1){25}} } \put(85,-55){\line(1,-1){25}}{\linethickness{0.5mm}\put(135,-55){\line(-1,-1){25}}}
 
 \put(135,-55){\line(1,-1){25}}  { \linethickness{0.5mm}{\put(185,-55){\line(-1,-1){25}}}} \put (188,-55){$\beta_3$}
 
       {\linethickness{0.5mm}   \put(60,-80){\line(0,-1){25}} \put(110,-80){\line(0,-1){25}}\put(160,-80){\line(0,-1){25}} 
              \put(60,-105){\line(-1,-1){15}}  \put(110,-105){\line(-1,-1){15}}  \put(160,-105){\line(-1,-1){15}} }
                      \put(60,-105){\line(1,-1){15}}     \put(110,-105){\line(1,-1){15}}  \put(160,-105){\line(1,-1){15}}

                           \put(207,-20){$\equiv$} 
             
                 \put(242,45){$\alpha_n$}            \put(255,45){\line(1,0){25}}
                  \put (230,15){$\alpha_{n-1}$}     \put(255,15){\line(1,0){25}}  \put(255,-25){$\vdots$} 
                  \put (242,-55){$\alpha_2$}  \put(255,-55){\line(1,0){25}}  
                  \put (242,-85){$\alpha_1$}  \put(255,-85){\line(1,0){25}}
                  
                  { \linethickness{1mm}{\put(280,52){\line(0,-1){145}}}}
                  
                    \put(280,-90){\line(1,-1){15}}\ \put(295,-110){$\xi_1$} 
                    
                       \put(280,30){\line(1,0){35}}\put(285,35){{\tiny $x_{n-1}^{(n-1)}$}}  
                           \put(280,0){\line(1,0){35}}\put(285,5){{\tiny $x_{n-2}^{(n-1)}$}} 
                            \put(295,-20){$\vdots$} 
                            \put(280,-40){\line(1,0){35}} \put(285,-35){{\tiny $x_{2}^{(n-1)}$}} 
                            \put(280,-70){\line(1,0){35}}\put(285,-65){{\tiny $x_{1}^{(n-1)}$}} 
                           
                          { \linethickness{1mm}{\put(315,37){\line(0,-1){115}}}}
                           \put(315,-75){\line(1,-1){15}}\ \put(330,-95){$\xi_2$}  

                          \put(315,15){\line(1,0){35}}\put(320,20){{\tiny $x_{n-2}^{(n-2)}$}} 
                           \put(315,-15){\line(1,0){35}}\put(320,-10){{\tiny $x_{n-3}^{(n-2)}$}} 
                            \put(330,-35){$\vdots$} 
                                             \put(315,-55){\line(1,0){35}}\put(320,-50){{\tiny $x_{1}^{(n-2)}$}} 
                            
                              { \linethickness{1mm}{\put(350,22){\line(0,-1){85}}}}  
                                \put(350,-60){\line(1,-1){15}}\ \put(365,-80){$\xi_3$}  
                                
                               \put(360,-5){$\cdots$}  \put(360,-40){$\cdots$} 
                               { \linethickness{1mm}{\put(385,7){\line(0,-1){55}}}}
                                \put(385,-45){\line(1,-1){15}}\ \put(400,-65){$\xi_{n-2}$} 
                               
                                \put(385,-5){\line(1,0){35}}\put(395,0){{\tiny $x_{2}^{(2)}$}}
                                  \put(385,-35){\line(1,0){35}}\put(395,-30){{\tiny $x_{1}^{(2)}$}}
                               
                                { \linethickness{1mm}{\put(420,0){\line(0,-1){40}}}}
                                   \put(420,-37){\line(1,-1){15}}\ \put(435,-57){$\xi_{n-1}$}  
                                
                                 \put(420,-20){\line(1,0){20}}\put(442,-22){{$\xi_n$}}             
 \end{picture}          \vskip45mm      
}
\caption{A Knutson--Tao's honeycomb  for $\gamma,\beta\gg \alpha, \xi=\gamma-\beta$ reduced to a  \GT triangle.
}
\label{HoneyComb-GTtriangle}
\end{figure}
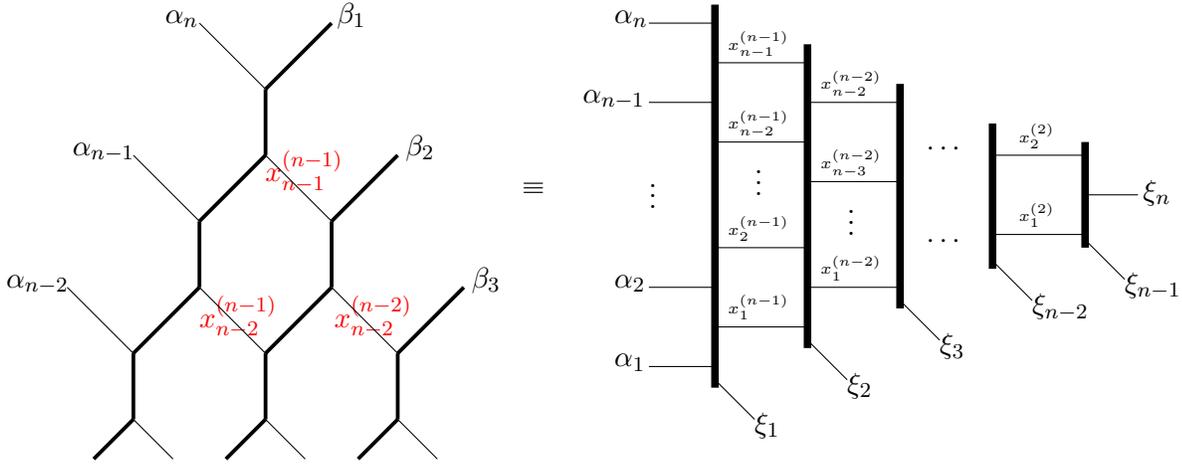

{\bigskip
Example. For $n=3$, these relations reduce  to 
\bea\nonumber \alpha_3\le x_2^{(2)} &\le& \alpha_2 \le x_1^{(2)} \le \alpha_1\\
 \label{poly3}  x_2^{(2)}\le \xi_3 &\equiv& x_1^{(1)} \le x_1^{(2)}\\
 \nonumber \sum_{i=1}^3\alpha_i  &=& \xi_1+x_1^{(2)}+x_2^{(2)} =\sum_{i=1}^3 \xi_i \,,
 \eea
 hence to the following bounds on, say, $ x_1^{(2)}$
 \be\label{intervalx22} \max(\alpha_2,\xi_2,\xi_3, \xi_2 {+} \xi_3-\alpha_2 )\le  x_1^{(2)}\le \alpha_1-\xi_1+\min(\xi_1, \alpha_2)\ee
 which leads to the expression (\ref{vol3}) of $\CI_{su(3)}(\alpha;\xi)$.}

\bigskip
A similar discussion may be carried out for the other algebras, based on the BZ inequalities. For example for $B_2$, see below sect.~\ref{B2}.


\subsection{The Kostka multiplicity as  a limit of LR coefficients}  
\label{KostkaFromLR}
In the same spirit as the large $\beta, \gamma$ limit above, we have 
\bea\mult_\lambda(\delta)= \lim_{s\to \infty} C_{\lambda\,\, \mu+s\rho}^{\mu+s\rho+\delta}\eea
independently of the choice of $\mu$. 
Actually we can choose $\mu = \rho$, so that the above formula reads
\bea \label{LRtoKostka} \mult_\lambda(\delta)= \lim_{s\to \infty} C_{\lambda\,\, s\rho}^{s\rho+\delta}\eea
 One can find an integer $s_c$ such that $\forall s\ge s_c, \,    \mult_\lambda(\delta)=  C_{\lambda\,\, s\rho}^{s\rho+\delta}$. 
 In the case $A_3$ for example,  taking
 $\lambda = \{4,5,3\}$, $\delta = \{-4,-2,5\}$,  one has $\mult_\lambda(\delta) = 26$, whereas taking $s = 4,5,6,7,8,9,\ldots$ one finds $C_{\lambda\,\, s\rho}^{s\rho+\delta} =  6, 19, 24, 26, 26, 26, \ldots$, so $s_c=7$.

 It is interesting to find a general upper bound for $s_c$ (we are not aware of any attempt of this kind in the literature). We shall find such a bound in the case of the Lie algebra $A_3$.
 As the reader will see, our method relies on a brute-force calculation; it would be nice to find a more elegant approach that could be generalized to $A_n$, or more generally, to any simple Lie algebra. 
For a given highest weight $\lambda$ and any weight $\delta$  belonging to its weight system, we calculate the Littlewood-Richardson coefficient $S(\lambda,\delta,s) =C_{\lambda\,\, s\rho}^{s\rho+\delta}$ from the Kostant-Steinberg formula, 
using the Kostant partition function for $A_3$ given in  sect.~\ref{kostantmult}, 
and show explicitly (this is indeed a brute-force calculation !) that the difference $S(\lambda,\delta,s+1) - S(\lambda,\delta,s)$ vanishes, equivalently, that $S(\lambda,\delta,s)$ is stationary, for values of $\lambda, \delta, s$ obeying the following set of constraints: 
$$ {\lambda_ 1} \geq 0, {\lambda_ 2} \geq 0, {\lambda_ 3} \geq 0,  {\lambda_ 1} + {\lambda_ 2} + {\lambda_ 3} > 0, s > 0,  4 (s + 1) > X_1 \geq 0 ,  2 (s + 1) > X_2 \geq 0,  4 (s + 1) > X_3 \geq 0.$$
where 
\begin{eqnarray*}
X_1&=& ( {\lambda_ 1} +2 {\lambda_ 2} + 3{\lambda_ 3}) - ({\delta_ 1} +2 {\delta_ 2} +3 {\delta_ 3}),\\
X_2&=& ( {\lambda_ 1} +2 {\lambda_ 2} + {\lambda_ 3}) - ({\delta_ 1} +2 {\delta_ 2} + {\delta_ 3}),\\
X_3 &=& ( 3{\lambda_ 1} + 2{\lambda_ 2} + {\lambda_ 3}) - (3{\delta_ 1} +2 {\delta_ 2} + {\delta_ 3}).
\end{eqnarray*}
The $s$-independent inequalities relating $\lambda$ and $\delta$ are nothing else than the Schur inequalities; equivalently, they can be obtained by writing that the partition $\alpha$ is larger than the weight $\xi$ for the dominance order on partitions --- as everywhere in this paper, $\alpha$ and $\xi$ refer to the Young components (partitions) associated with the weights $\lambda$ and $\delta$. 
In the previous set of constraints, one can replace the $s$-independent inequalities relating $\lambda$ and $\delta$ by the following ones:
$-(\lambda_1+\lambda_2+\lambda_3)\leq \phi \leq (\lambda_1+\lambda_2+\lambda_3)$, where $\phi$ can be $\delta_1$, $\delta_2$, $\delta_3$ or  $(\delta_1+\delta_2 +\delta_3)$. 
The latter, namely the one with $\phi = (\delta_1+\delta_2 +\delta_3)$, expresses the fact that $\lambda$ (resp. $ - \overline {\lambda}$) is the highest weight (resp. is the lowest weight).
The obtained $s$-dependent inequalities imply  $4(s_c+1) < \max(X_1, 2 X_2, X_3)+1$.

For a given highest weight $\lambda$ this  allows us to obtain a bound independent of the choice of the weight $\delta$. One finds $s_c \leq 2(\lambda_1+\lambda_2+\lambda_3)$.
Indeed, one can check explicitly that, given $\lambda$,  and for any $\delta$ of its weight system, the function $S(\lambda,\delta,s)$ is stationary for $s \geq 2\sum_j \lambda_j$.


\subsection{Reduced O-blades and reduced isometric honeycombs}   
\label{ObladesVersusLianas}

The general discussion carried out in section \ref{polytope} could be expressed in terms of other pictographs, for instance BZ-triangles, O-blades, or isometric honeycombs (see our discussion in \cite{CZ0} or \cite{CZ1} for a presentation of the last two).
We have seen how \KT honeycombs are ``reduced'' to \GT patterns when one moves from the Horn problem to the Schur problem. An analogous reduction holds if we use isometric honeycombs or rather, their O-blades partners. 
We shall illustrate this with SU(4) by choosing the (dominant) weight $\lambda = \{4,5,3\}$ and the weight $\delta = \{-4,-2,5\}$;
here components are expressed in terms of Dynkin labels\footnote{Remember that external sides of \KT honeycombs are labelled by integer partitions whereas external sides of O-blades or of isometric honeycombs are labelled by the Dynkin labels of the chosen weights.}.
Equivalently, in terms of partitions, we have $\alpha = (12,8,3,0)$, and $\xi = (-1, 3, 5, 0) + (4,4,4,4) = (3, 7, 9, 4)$.

$\lambda$ is the highest weight of an irreducible representation (of dimension $16500$), and the weight subspace associated with $\delta$ has dimension $26$. 
There is a basis in this representation space for which every basis vector can be attached to a semi-standard Young tableau with filling $1,2,3,4$ and shape (partition)  $\alpha$,   
or, equivalently, to a Gelfand--Tsetlin pattern. Here is one of them, and one  sees immediately that its associated basis vector indeed belongs to the weight subspace defined by $\delta$ (or $\xi$):

{\centering
Young tableau : 
$
{\footnotesize
\begin{array}{cccccccccccc}
 1 & 1 & 1 & 2 & 2 & 2 & 2 & 3 & 3 & 3 & 3 & 3 \\
 2 & 2 & 2 & 3 & 3 & 3 & 4 & 4 & \text{} & \text{} & \text{} & \text{} \\
 3 & 4 & 4 & \text{} & \text{} & \text{} & \text{} & \text{} & \text{} & \text{} &
   \text{} & \text{} \\
\end{array}
}
$
\hspace{0.5cm}
with \GT pattern :
\hspace{0.5cm}
$
{\footnotesize
\begin{array}{c}
12,8,3,0 \\
12,6,1 \\
7,3 \\
3 \\
\end{array}
}
$
\par}

Remember that the sequence of lines of the \GT pattern is obtained from the chosen Young tableau by listing the shapes of the tableaux obtained by removing successive entries, starting from the largest one (here $4$):
{\footnotesize  $\{\{1,1,1,2,2,2,2,3,3,3,3,3\},\{2,2,2,3,3,3\},\{3\}\}$, $\{\{1,1,1,2,2,2,2\},\{2,2,2\}\}$,  $\{\{1,1,1\}\}$}. 

Using $s = 1000$ in (\ref{LRtoKostka}), we display in (Fig. \ref{ObladeSchurExample}, Left) one\footnote{As we shall see below, this choice (among the $26$ O-blades obtained when $s>s_c$) is not arbitrary but dictated by our wish to establish a link with the previously chosen Young tableau.} of the O-blades describing the space of intertwiners for the triple\footnote{In order to ease the discussion of the  correspondence with Young tableaux,  and also to draw the weight $\lambda$ on the bottom of each pictograph, it is better to display  the O-blades (or their isometric honeycombs partners) associated with the Littlewood-Richardson coefficient  $C_{s \rho + \delta, s \rho}^\lambda$  or with $C_{s \rho, s \rho + \delta}^\lambda$
rather than those associated with $C_{\lambda, s \rho}^ {s \rho + \delta}$. Notice that the numbers of such pictographs are all equal to $\mult_\lambda(\delta)$ when $s$ is big enough. The choice made in Fig. \ref{ObladeSchurExample}, \ref{allthereducedObladesNew} corresponds to the triple $({(s \rho, s \rho + \delta)} \mapsto \lambda)$.} 
$({(s \rho, s \rho + \delta)} \rightarrow \lambda))$. 
It is a member of a one-parameter ($s$) family of O-blades for which the integers carried by the horizontal edges are $s$-dependent and for which the other edge values stay constant as soon as $s>s_c$. What remains in the limit of large $s$ (meaning $s>s_c$) is only the ``reduced O-blade'' given in  (Fig. \ref{ObladeSchurExample}, Middle) where we have removed the 0's and the values carried by the $s$-dependent {\sl horizontal} edges since they are irrelevant. 
\\
The net result is that we have as many distinct reduced O-blades (for instance those obtained by reducing the ones associated  with the  triple ${((s\rho, s\rho + \delta)} \rightarrow \lambda$)) as  we have distinct \GT patterns, namely,  a number equal to the dimension $\mult_\lambda(\delta)$ of the weight subspace defined by $\delta$.  Those ($26$ of them) associated with the example chosen above are displayed on Fig.~\ref{allthereducedObladesNew}.
\\
 Isometric honeycombs cannot be displayed on a page if $s$ is large (because they are isometric !), so we choose $s=8$ and display on Fig.~\ref{ObladeSchurExample} (Right)  the isometric honeycomb partner of the O-blade obtained, in the same one-parameter family, for this value of $s$ (here $s_c = 7$ so that taking $s=8$ is enough). 
Now the integers carried by the {\it vertical} edges are irrelevant (they would change with $s$): what matters are the integers carried by the non-vertical edges of the --possibly degenerate-- parallelo-hexagons.

 The particular reduced O-blade displayed in Fig.~\ref{ObladeSchurExample} (Middle) can also be obtained directly,  by a simple combinatorial rule, from the Young tableau displayed previously, or from the corresponding \GT pattern, 
without any appeal to the limit procedure (\ref{LRtoKostka}) using Littlewood-Richardson coefficients, 
and without considering their associated pictographs:  see the Appendix ``Lianas and forests''.

 \begin{figure}[htb]
  \centering
       \includegraphics[width=20pc]{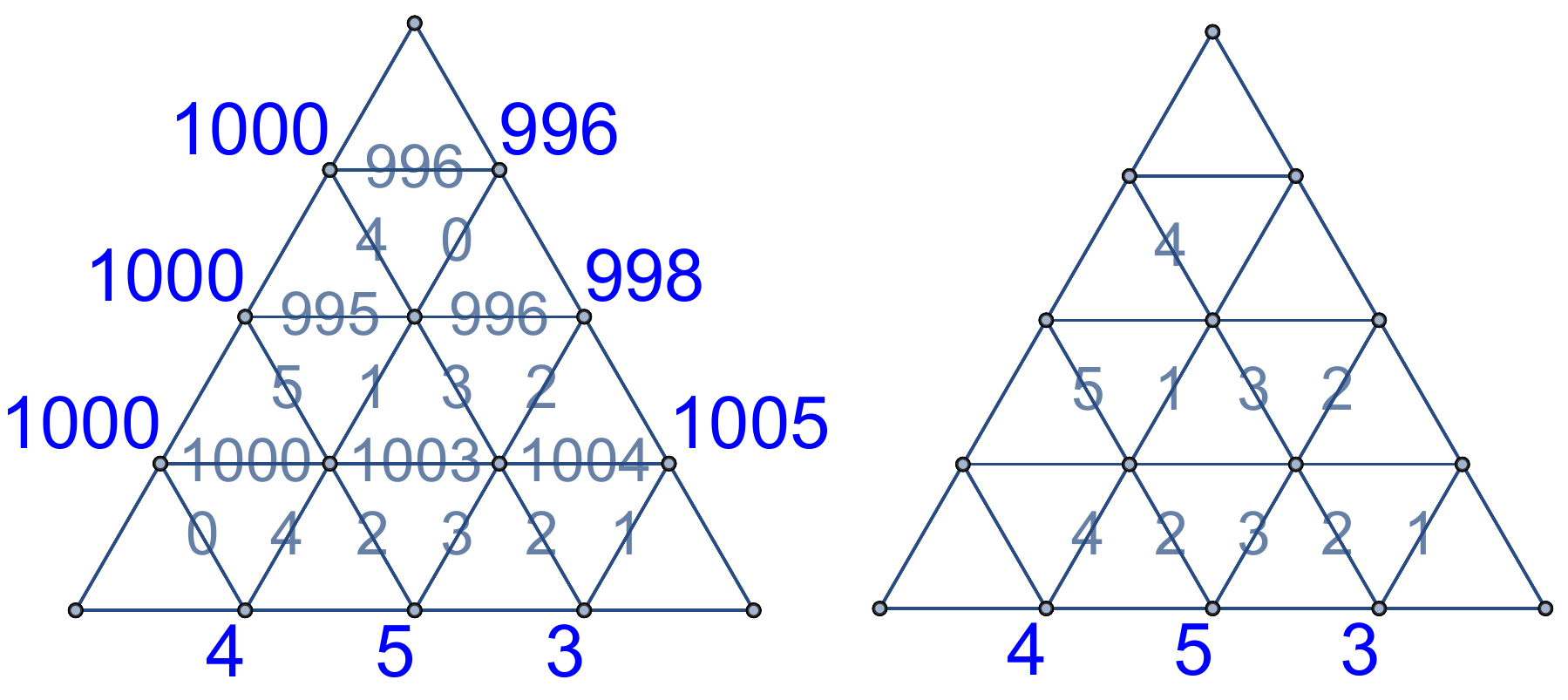}
       \includegraphics[width=7pc]{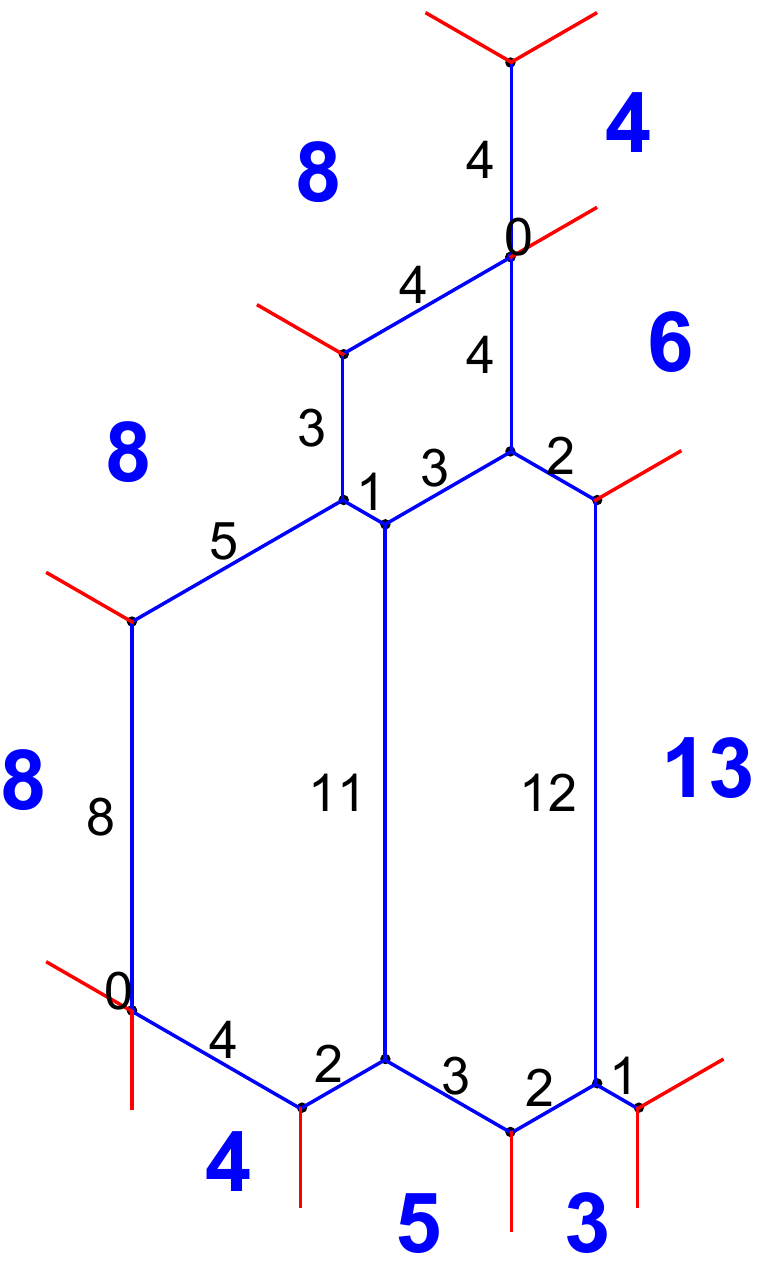}
       \caption{
       Left : One of the $26$ O-blades associated with  $C_{s \rho + \delta, s \rho}^\lambda$ for $\lambda = \{4,5,3\}$, $\delta=\{-4,-2,5\}$, $s = 1000$. 
       Middle: same, with the $0$'s and the $s$-dependent labeling of horizontal lines removed. 
       Right: The corresponding isometric honeycomb for the choice $s=8 > s_c$.} 
        \label{ObladeSchurExample}
       \end{figure}


\subsection{Asymptotics of the $\CI$--multiplicity formula:  $\CI$ as a volume and stretching polynomials.}
\label{Asymptotics}
We may repeat the same chain of arguments as in Horn'problem \cite{CZ1, CMSZ1}: 
upon scaling  by $p\gg 1$ of $\lambda, \delta$  in (\ref{I-mult})\\
\bea\nonumber  \CI(p\lambda+\rho,p \delta)&\approx& \CI(p\lambda,p \delta)=p^{(n-1)(n-2)/2} \CI(\lambda,\delta)\\
\nonumber &=& \sum_{\kappa\in K} r_\kappa \sum_\tau C_{p\lambda\, \kappa}^\tau \mult_\tau(p\delta)\\
\nonumber &\approx& \sum_{\kappa\in K} r_\kappa \sum_\tau \sum_{k\in[\kappa]}\delta_{\tau,p\lambda+k} \mult_\tau(p\delta)\\
\nonumber &=& \sum_{\kappa\in K} r_\kappa  \sum_{k\in[\kappa]} \mult_{p\lambda+k}(p\delta)\\
\nonumber &=& \underbrace{\sum_{\kappa\in K} r_\kappa  \dim_\kappa}_{=1}  \mult_{p\lambda}(p\delta)\\
 &=& \mult_{p\lambda} (p\delta)=  p^d \mathrm{vol}_d(\CP(\lambda;\delta))+\cdots
 \label{asy-scaling}
 \eea
where  $\CP(\lambda;\delta)$  is the polytope defined in sect.~\ref{polytope} or associated with one of the pictographs mentioned 
in sect.~\ref{ObladesVersusLianas}. Thus for generic cases for which $d=(n-1)(n-2)/2$, we have the identification
 \be  \CI(\lambda,\delta) = \mathrm{vol}_d(\CP(\lambda;\delta))\,, \ee
 (while for the non generic cases, both $\CI$ and  $\mathrm{vol}_{(n-1)(n-2)/2}$ vanish).
 
 \bigskip
 For illustration we consider the weights $\lambda=\{4,5,3\}$ {\sl and\/} $\delta = \{-4,-2,5\}$ of SU(4) already chosen in a previous section.
The multiplicities obtained by scaling them by a common factor $p=1,2,3,\ldots$ lead to the sequence  $26, 120, 329, 699, 1276, 2106, 3235, 4709, 6574, 8876, \ldots$, which can be encoded by the cubic polynomial 
$\tfrac{1}{6} (6 + 35 p + 69 p^2 + 46 p^3)$ whose dominant term is $23/3$,
which is indeed the value of the Schur volume function $\CI(\lambda, \delta)$.

Another way to obtain this value is to use the volume function {\sl for the Horn problem}, that we called ${\mathcal J}(\lambda, \mu, \nu)$  in refs \cite{CZ1} and \cite{CMSZ1}.  In the generic case this function gives the volume of the hive polytope associated with 
the triple $(\lambda, \mu)\rightarrow \nu$ \ie the polytope whose $C_{\lambda, \mu}^\nu$ integer points label the honeycombs for this specific space of intertwiners. It also gives the leading coefficient of the Littlewood-Richardson polynomial 
(a polynomial\footnote{We assume here that the chosen Lie group is SU(n), otherwise, this object may be a quasi-polynomial.} in the variable $p$) giving, when $p$ is a non-negative integer, the LR coefficient $C_{p \lambda, p\mu}^{p\nu}$ for highest weights scaled by $p$.
 Since the Kostka numbers (multiplicities of weights) can be obtained as a limit of LR coefficients for special arguments (see~(\ref{LRtoKostka})), the same is true under scaling.
For $s>s_c$ (determined by the choice of $\lambda$ and $\delta$), we have therefore $\CI(\lambda, \delta) = {\mathcal J}(\lambda, s \rho, s \rho + \delta)$. We check, on the same example as before, that we have indeed, 
${\mathcal J}(\{4,5,3\},\{s,s,s\},\{s-4,s-2,s+5\}) = 13/3, 7, 23/3, 23/3, 23/3, 23/3, \ldots $ for $s=5, 6, 7, 8, \ldots$ As expected,  the values of ${\mathcal J}$ stabilize, and the asymptotics, \ie $\CI(\lambda, \delta)$,  is reached for $s =  7$.
Notice the various kinds of scalings involved here: 1) a scaling of $\rho$ by $s$, with the weights $\lambda$ and $\delta$ remaining constant,  2) a scaling of $\lambda$ and $\delta$ by the non-negative integer $p$, 
3) a simultaneous scaling of the three arguments of ${\mathcal J}$  by a non-negative integer $p$ giving rise to a polynomial (in $p$) encoding  the LR coefficients $C_{p \lambda, p s \rho}^{p s\rho + p\delta}$.

\section{The case of $B_2$}
\label{B2}

\begin{figure}[htb]
\begin{center}
\includegraphics[width=0.3\textwidth]{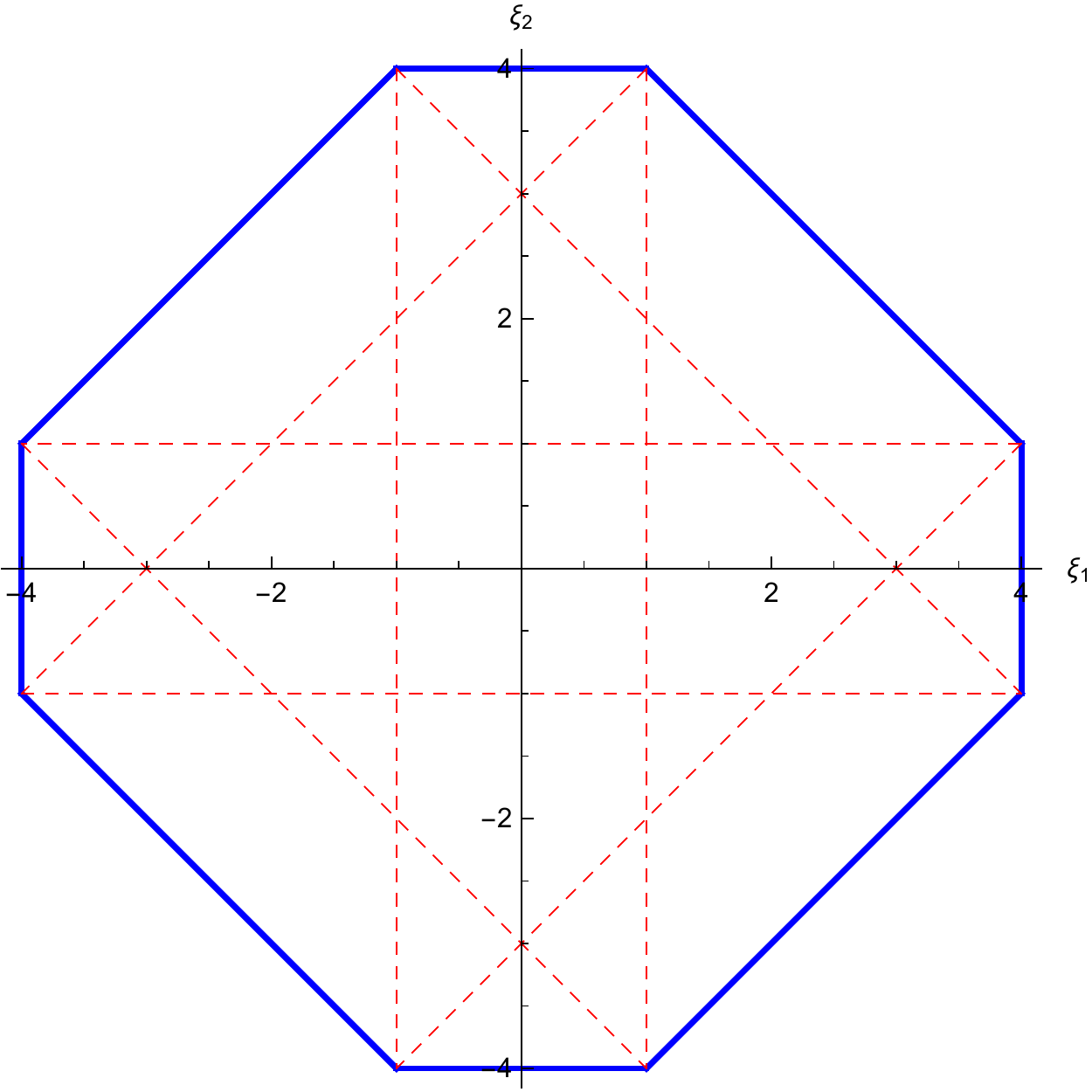}\qquad \includegraphics[width=0.3\textwidth]{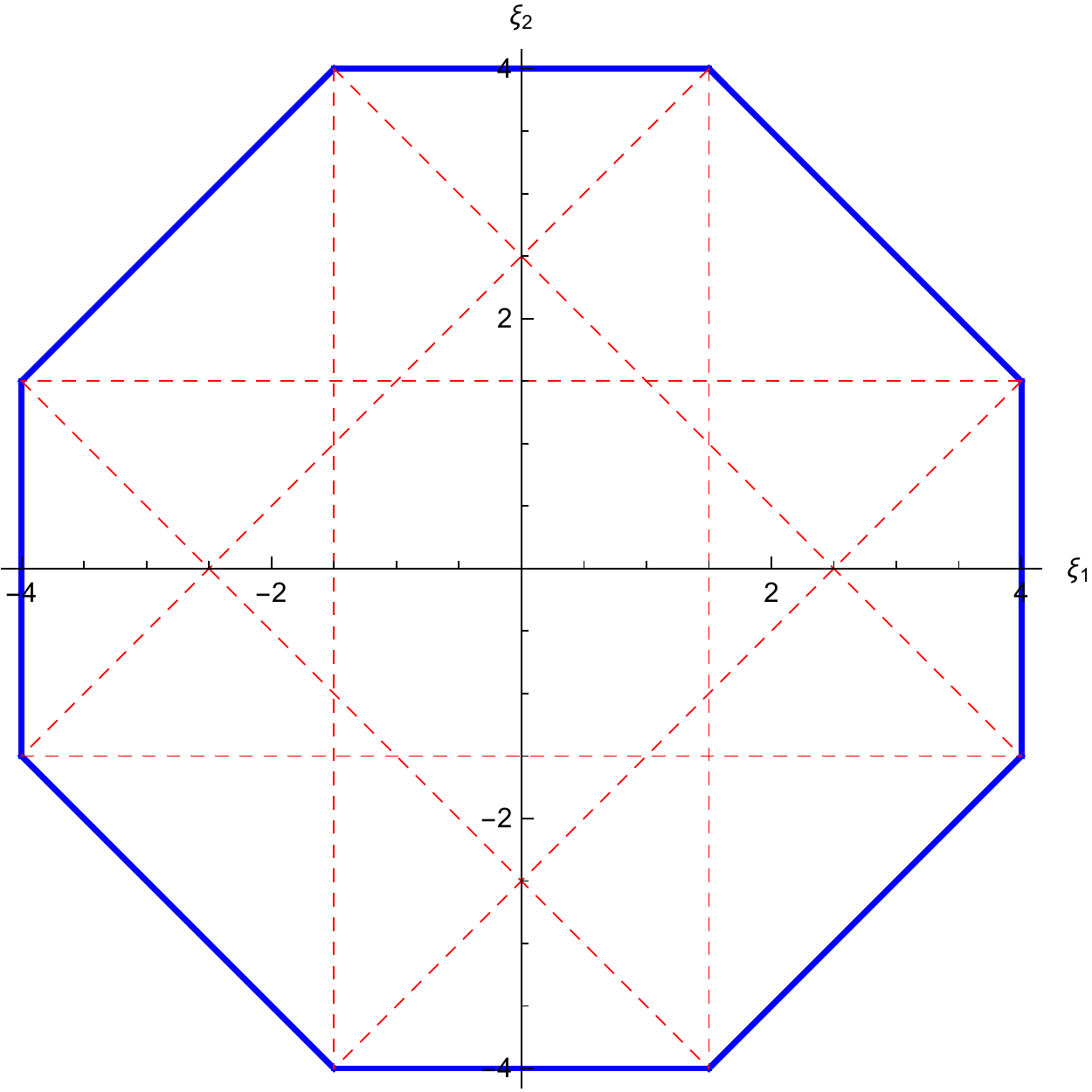}\qquad \includegraphics[width=0.3\textwidth]{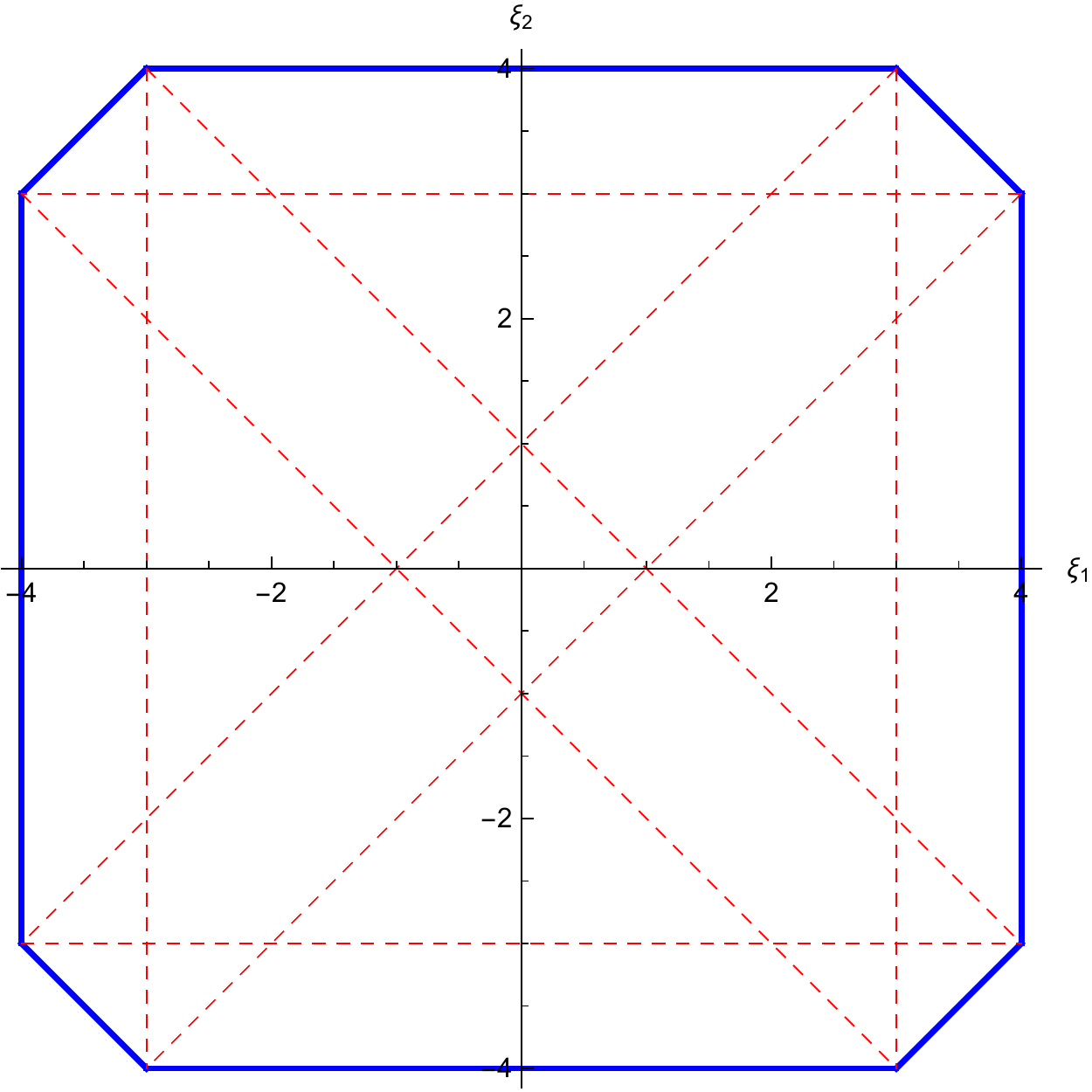}
\caption{The Schur polytope $\mathcal{O}$ in the $B_2$ case, with its singular lines, in the three cases 
$\alpha=(4,1)$, $\alpha=(4,\frac{3}{2})$ and  $\alpha=(4,3)$.  }
\label{Schur-PolyB2}  
\end{center}
\end{figure}

\subsection{The Schur octagon}
 Consider the orbit $\CO_\alpha$ of the group $\SO(5)$ acting by conjugation on a block-diagonal skew symmetric matrix \\
$A=\diag\bigg( {\footnotesize\begin{bmatrix}0&\alpha_i\\ -\alpha_i&0 \end{bmatrix}}_{i=1,2}, 0\bigg)$ with real $\alpha_i$. 
Note that one may choose $0\le \alpha_2\le \alpha_1$. Schur's problem reads: 
what can be said about the projections  $\xi_i$ of a matrix of $\CO_\alpha$ on an orthonormal  basis $X_i$ of the Cartan algebra,
$\xi_i=\tr O.A.O^T.X_i$, $i=1,2$, where $O \in \SO(5)$? 
More specifically, if the matrix $O$ is taken randomly and uniformly distributed in SO(5), what is the PDF of the $\xi$'s? 
There is again a connection with the corresponding Horn problem, as discussed for example in \cite{CMSZ1}.
Taking the limit $\beta,\gamma\gg \alpha$, $\xi=\gamma-\beta$ finite, in the expressions and figures of
\cite{CMSZ1}, we find that the support of the Schur volume is the octagon
\be\label{B2-octagon}   
\mathcal{O}\quad:\qquad-\alpha_1\le \xi_1, \xi_2 \le \alpha_1\,, \qquad 
-(\alpha_1+\alpha_2)\le \xi_1\pm \xi_2 \le \alpha_1+\alpha_2\,, \ee
and that the singular lines of the PDF, or of the associated volume function,
are $$\xi_1,\xi_2=\pm \alpha_2\,, \qquad\xi_1\pm \xi_2= \pm(\alpha_1-\alpha_2)\,,$$
see Fig. \ref{Schur-PolyB2}. Three cases occur, depending on whether the ratio $\alpha_2/\alpha_1$ belongs to
$(0,1/3),\, (1/3, 1/2)$ or $(1/2,1)$. For a given value of that ratio, $1+4 \times 6=25$ cells occur, and in total, 
one should have $1+4 \times 8=33$ possible cells, in accordance with Bliem's result~\cite{Bl}. 

Working out this same limit in the expressions of the Horn volume given in \cite{CMSZ1}, we find a (relatively) simple 
expression for the Schur volume function
\bea \nonumber \CI_{B_2}(\alpha;\xi ) &=& 0 \qquad \mathrm{if}\  \xi \notin \mathcal{O} \\
\label{DeltaIB2} \Delta\CI_{B_2}(\alpha;\xi)   &=& \begin{cases}  \pm \oh(\xi_i-\alpha_j)^2 & \mathrm{across\ any\ line} \ \xi_i=\alpha_j \qquad i,j=1,2\\
\pm \inv{4}(\xi_1+\epsilon \xi_2 -\epsilon''(\alpha_1+\epsilon' \alpha_2))^2 & \mathrm{across\ any\ line}\ \xi_1+\epsilon \xi_2 =\epsilon''(\alpha_1+\epsilon'\alpha_2)
\end{cases}
\eea
$\epsilon, \epsilon',\epsilon'' =\pm 1$, and 
with the overall sign of the change determined by the prescription of Fig. \ref{four-prong}. 

  \begin{figure}[htb]
  \centering
       \includegraphics[width=0.15\textwidth]{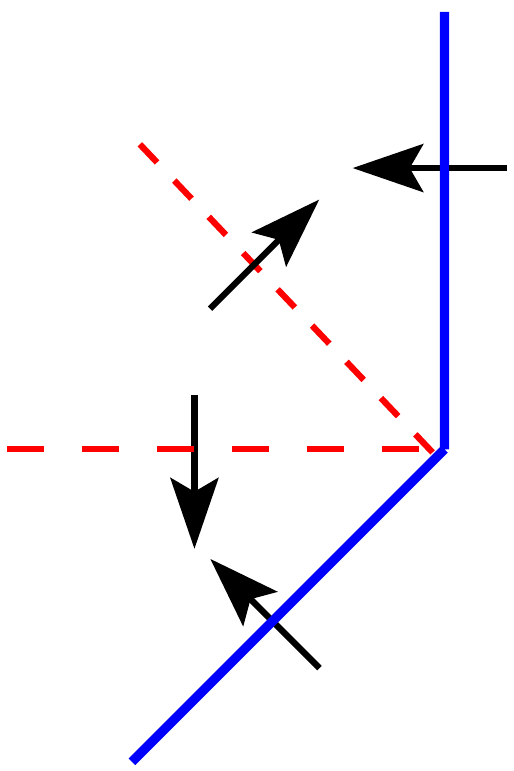}  \caption{\label{four-prong} Prescriptions for the changes of $\CI_{B_2}$ across the lines
       emanating from a  vertex: the sign in (\ref{DeltaIB2}) is $+$ if the line is crossed along the direction of the arrow, $-$ otherwise. 
       This holds for any rotated configuration of that type. }  \end{figure}

One may finally plot the resulting PDF $\p(\xi|\alpha)=\frac{3}{2\alpha_1\alpha_2(\alpha_1^2-\alpha_2^2)} \CI_{B_2}(\alpha,\xi)$ for 
a given value of $\alpha$ and compare it with the ``experimental" histogram of a large random sampling of matrices of $\CO_\alpha$, 
see Fig. \ref{histo-B2}. }

      \begin{figure}[tbh]
  \centering \includegraphics[width=0.4\textwidth]{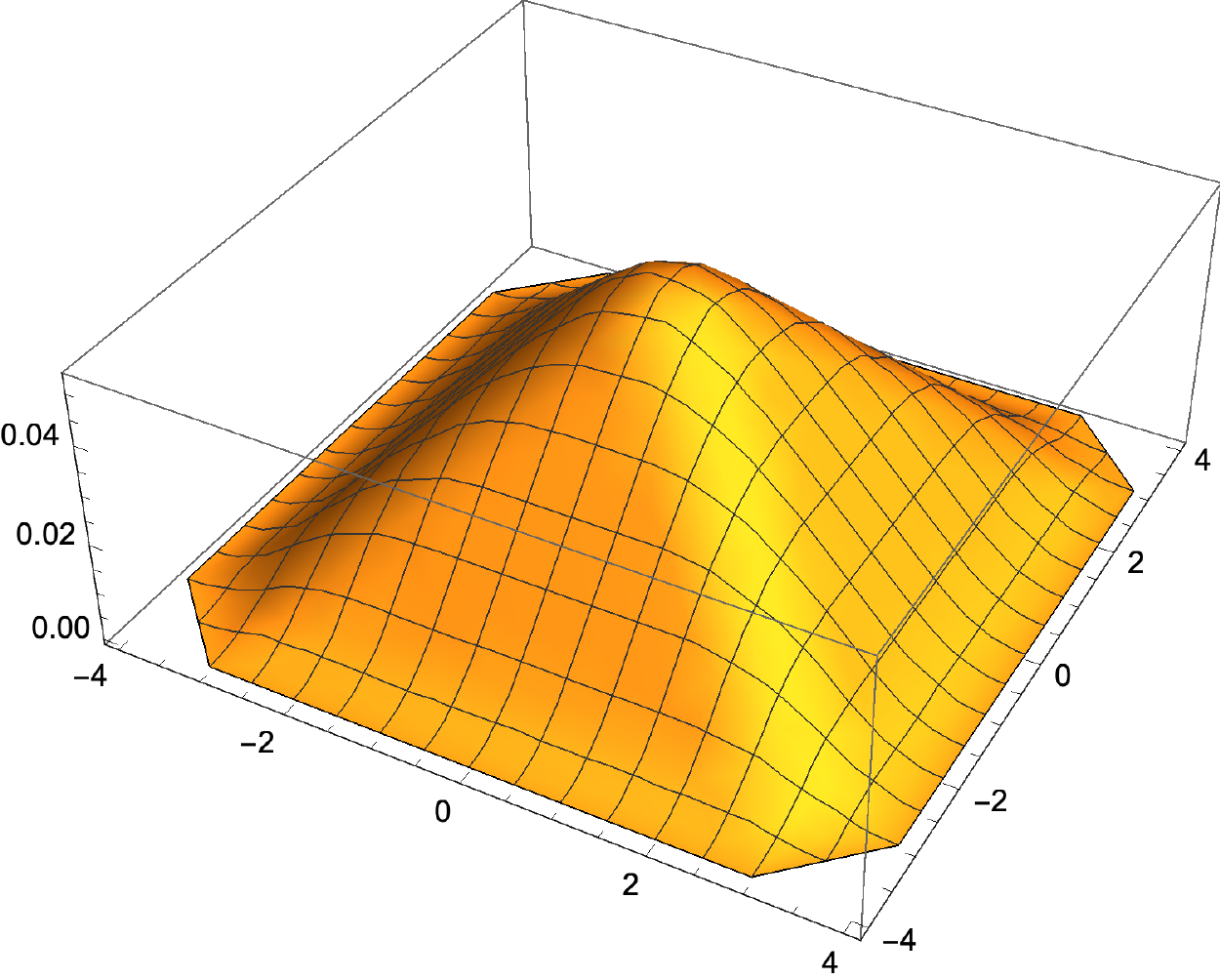}  \qquad
       \includegraphics[width=0.4\textwidth]{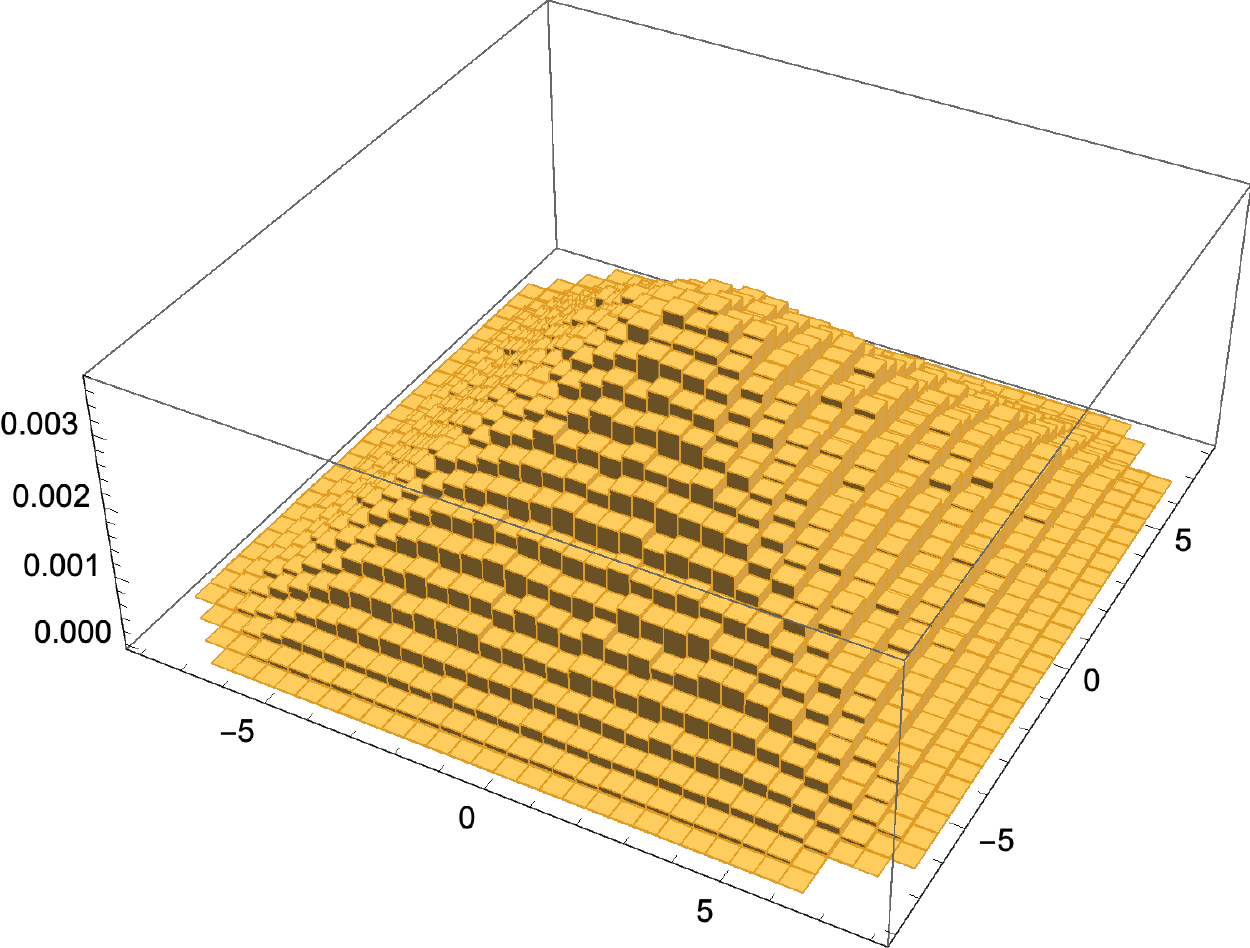}  \caption{\label{histo-B2} 
       Comparing the plot of the PDF with the histogram of $\xi$ values obtained from $10^6$ matrices of the orbit
       $\CO_\alpha$, for $\alpha=(4,3)$. }  \end{figure}

Repeating the steps followed in sect. \ref{CI-mult} and making use of results in sect. 3 of \cite{CMSZ1}, 
one derives a relation between $\CI_{B2}$ and multiplicities: 
\be\label{CI-mult-B2} \CI_{B_2}(\lambda';\delta) =\begin{cases}\inv{8} \( 3\, \mult_\lambda(\delta) +\sum_\tau  C_{\lambda\,\{1,0\}}^\tau\mult_{\tau}(\delta)  \)& \mathrm{if}\ \lambda-\delta \in Q \cr
\inv{4} \sum_\tau C_{\lambda\,\{0,1\}}^\tau\mult_{\tau}(\delta)& \mathrm{if}\ \lambda-\delta-\rho \in Q \,.
\end{cases}
\ee
In particular for $\lambda=0$, 
\be \CI_{B_2}(\rho;\delta) =\begin{cases}\inv{8} \( 3\, \delta_{\delta 0} +\mult_{\{1,0\}}(\delta)  \)& \mathrm{if}\ \delta \in Q \cr
\inv{4}\mult_{\{0,1\}}(\delta)  & \mathrm{if}\ \delta-\rho \in Q \
\end{cases}
\label{IB2rhodelta}
\ee
For example, take for $\delta$ the short simple root, $\delta=
(0,1)=\{-1,2\}\in Q$,  $\CI_{B_2}(\rho;\delta)=\inv{8}$, $\mult_{\{1,0\}}(\delta) =1$.

  {\bf Remarks}. \\
   1. A general piece-wise polynomial expression of $\mult_\lambda(\delta)$ has been given by Bliem \cite{Bl}. His expressions 
   should be consistent with  the relation  (\ref{CI-mult-B2}) and our implicit expression of $\CI_{B_2}$ in (\ref{DeltaIB2}). \\
2. One should also notice that (\ref{LRtoKostka}) holds true in general, and thus in the current $B_2$ case. For example, 
 $\mult_{\{20, 12\}} (\{18, -6\})=56$ may be  recovered asymptotically as 
$C_{\{20,12)\}\,\{s,s\}}^{\{s+18, s-6\}}$, which takes the values 
${0,0,0,0,0,3,8,14,20,26,31,36,40,44,47,50,52,54,55,56,56,56,56,56,56}$ 
as $s$ grows from 1 to 25. Hence here $s_c = 20$.


 \subsection{$A_3$ versus $B_2$}
It is well known that the $B_2$ 
root system may be obtained   by folding that of $A_3$. It is thus suggested to compare Kostka multiplicities for cases that enjoy some 
symmetry in that folding. Consider in particular the special case of $\lambda=\{\lambda_1,1,\lambda_1\} $.
       Inequalities on the three parameters $i,j,k$ of $A_3$ oblades for $\mathrm{mult}_{\{\lambda_1,\lambda_2=1,\lambda_1\}}(\{\delta_1,\delta_2,\delta_1\})$ reduce to
\bea && \{i\leq {\lambda_1},\  {\delta_2}\leq 2 k+1,\  {\delta_1}+k\leq i,\  {\delta_1}+j+k+{\lambda_1}\geq
   i,\   -1\leq {\delta_2}-2 (j+k)\leq 1, 
   \\ \nonumber &&
   i\geq 0,\  i+j\geq 0,\  {\delta_1}+j+k\leq i,\ 
   {\delta_1}+j+2 k\leq i,\  {\delta_2}+1\geq 2 j,\  i\geq k\}\eea
Consideration of a large number of examples then suggest the following\\
{\bf Conjecture 1}:\\
{\sl i)  the number of triplets $(i,j,k)$  satisfying these inequalities is a square integer, viz 
$$\mathrm{mult}_{\{\lambda_1,1,\lambda_1\}}(\{\delta_1,\delta_2,\delta_1\}) =m^2\,;$$
ii) the corresponding $B_2$ multiplicity is then 
$\mathrm{mult}_{\{\lambda_1,1\}}(\{\delta_1,\delta_2\}) =m$\,.}

\subsection{Stretching (quasi)polynomial for the $B_2$ multiplicity}

Let $\kappa=\lambda-\delta$. As well known, $\mult_\lambda(\delta)$ vanishes if $\kappa \notin Q$, the root lattice, \ie
if $\kappa_2$ (in Dynkin indices) is odd. 
Now consider the stretched multiplicity $\mult_{s\lambda}(s\delta)$, $s\in \N$. It is in general a quasi-polynomial of $s$. 
Again, we have found a fairly  strong evidence for, and we propose the following \\
{\bf Conjecture 2}:\\
{\sl i)  for $\lambda_2$ and $\delta_2$ both even, 
   $\mult_{s\lambda}(s\delta)$ is a polynomial of $s$ for $\kappa_1$ even, and,
   except for a finite number of cases, a quasi-polynomial for $\kappa_1$ odd;\\ 
 ii) for $\lambda_2$ and $\delta_2$ both odd, 
generically it is a quasi-polynomial, except if  $2 \CI_{B_2}$  (twice the Schur volume) is an integer.}\\

 \subsection{Pictographs for $B_2$ ?}
 A combinatorial algorithm based on Littelmann's paths~\cite{Li} has been proposed by Bliem\cite{Bl}.
{We tried (hard) to invent a $B_2$ analog of O-blades, either degenerate (Kostka coefficients) or not (LR coefficients).
 It seems that it cannot be done without introducing edges carrying both positive and negative integers, and the result is not particularly appealing, so we  leave this as an open problem~!}


\section*{Acknowledgements}  It is a great pleasure to thank Colin McSwiggen for several enlightening discussions and suggestions.


\section*{Appendix :  Lianas and forests}

We now explain how to obtained directly a reduced O-blade from an $\SU(n)$ Young tableau.
We do not assume that the given tableau is semi-standard, but its entries along any chosen column should increase when going down.
\\
Draw an equilateral triangle with sides of length $n$, choose one side (``ground level'') and mark the inner points $1,\ldots, n-1$.
To every column (with $j$ elements)  of the chosen tableau, associate a zigzag line (a {\sl liana}) going upward, but only north-west or north-east, from the marked point~$j$ located on the ground level of the triangle. 
 When going up, there are $n+1$ consecutive levels for $\SU(n)$, the ground level being the first.
An entry marked $p$ in the chosen column makes the associated liana to grow upward to the left between levels $(n+1)-p$ and $(n+1)-(p-1)$, otherwise, it grows upward to the right.
When it reaches the boundary of the surrounding triangle, the liana continues upward, following the boundary.
One finally superimposes the lianas rooted in the same points (the positions $1,2,\ldots, n-1$), formally adding the characteristic functions of their sets, when they grow from the same point, and obtain in this way a {\sl liana forest}.\\
\noindent
As in  sect.~\ref{ObladesVersusLianas} we consider the $\SU(4)$  example
$
{\footnotesize
\begin{array}{cccccccccccc}
 1 & 1 & 1 & 2 & 2 & 2 & 2 & 3 & 3 & 3 & 3 & 3 \\
 2 & 2 & 2 & 3 & 3 & 3 & 4 & 4 & \text{} & \text{} & \text{} & \text{} \\
 3 & 4 & 4 & \text{} & \text{} & \text{} & \text{} & \text{} & \text{} & \text{} &
   \text{} & \text{} \\
\end{array}
}
$.
We have $12$ lianas with five levels. Reading this Young tableau from right to left we see that left-going directions occur between levels specified by
 {\footnotesize $\{\{3\},\{3\},\{3\},\{3\},\{3,4\},\{2,4\},\{2,3\},\{2,3\},\{2,3\},\{1,2,4\},\{1,2,4\},\{1,2,3\}\}$}, to which we associate the lianas displayed in Fig.~\ref{lianas}. 
 For instance the $6^{th}$ liana, which is rooted in the second point, is described by $\{2,4\}$, so it grows upward to the left above levels $3=5-2$ and $1=5-4$, otherwise it grows to the right.

 \begin{figure}[htb]
 \centering
 \includegraphics[width=30pc]{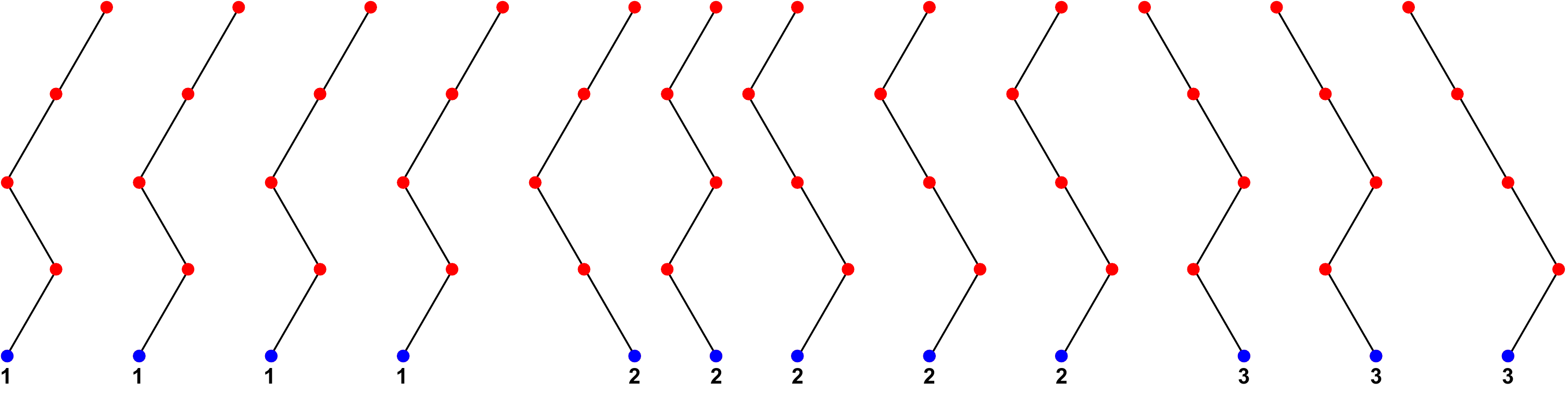}
\caption{SU(4):  Lianas associated with the Young tableau chosen in the text.} 
\label{lianas}
 \end{figure}

When superimposing the lianas rooted in the same points we recover the reduced O-blade displayed of Fig.~\ref{ObladeSchurExample}, (Middle), as a liana forest: see Fig.~\ref{lianaForest}.
Notice that the labels $9$, $3$, $5$, $3$, $1$, on the boundary of the triangle given in Fig.~\ref{lianaForest} (Right), are absent in the reduced O-blade,
but they can be obtained immediately from the latter by using the property that when a liana reaches the boundary, it continues upward, following the boundary.

 \begin{figure}[htb]
 \centering
  \includegraphics[width=12pc]{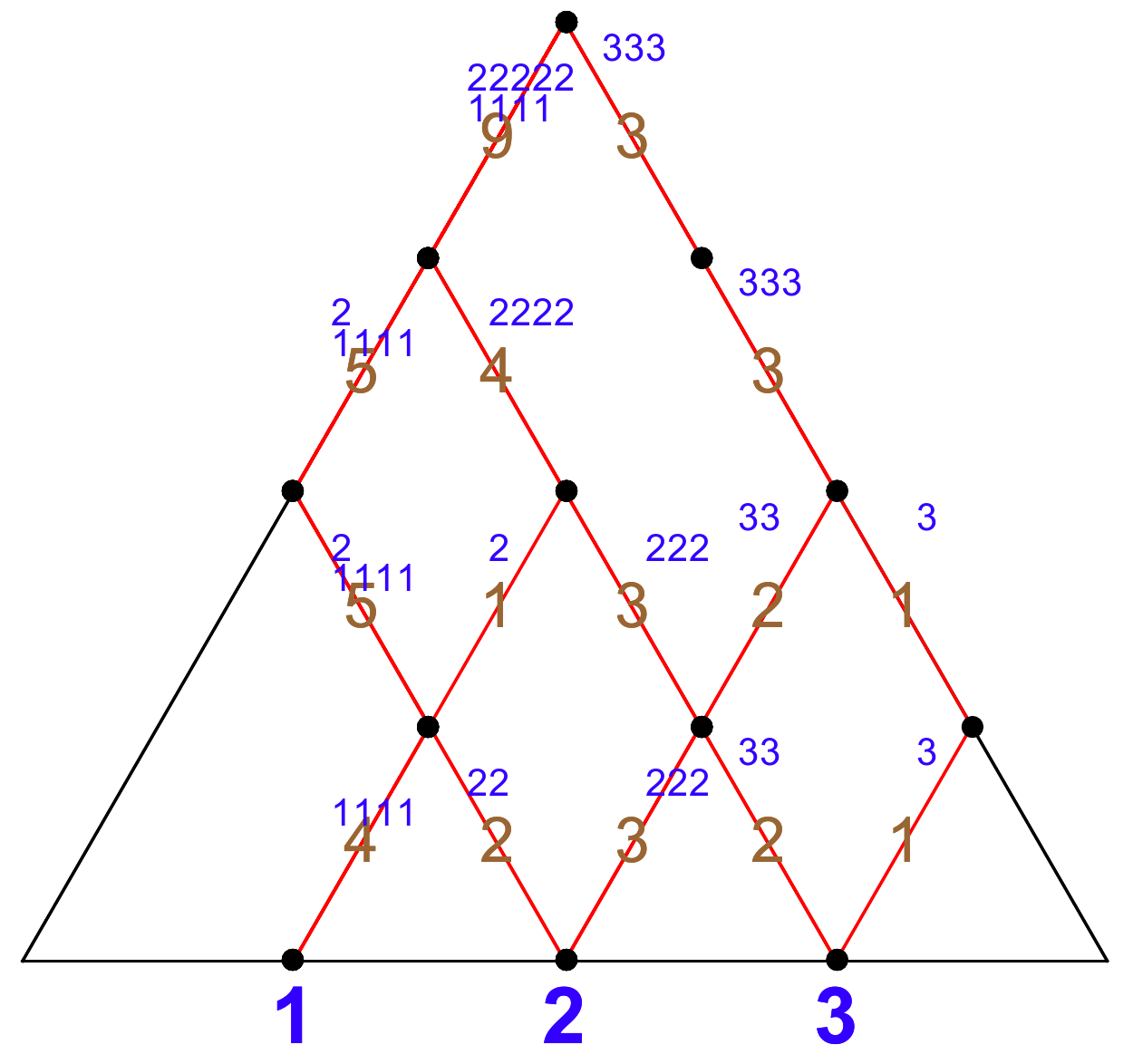}
    \includegraphics[width=12pc]{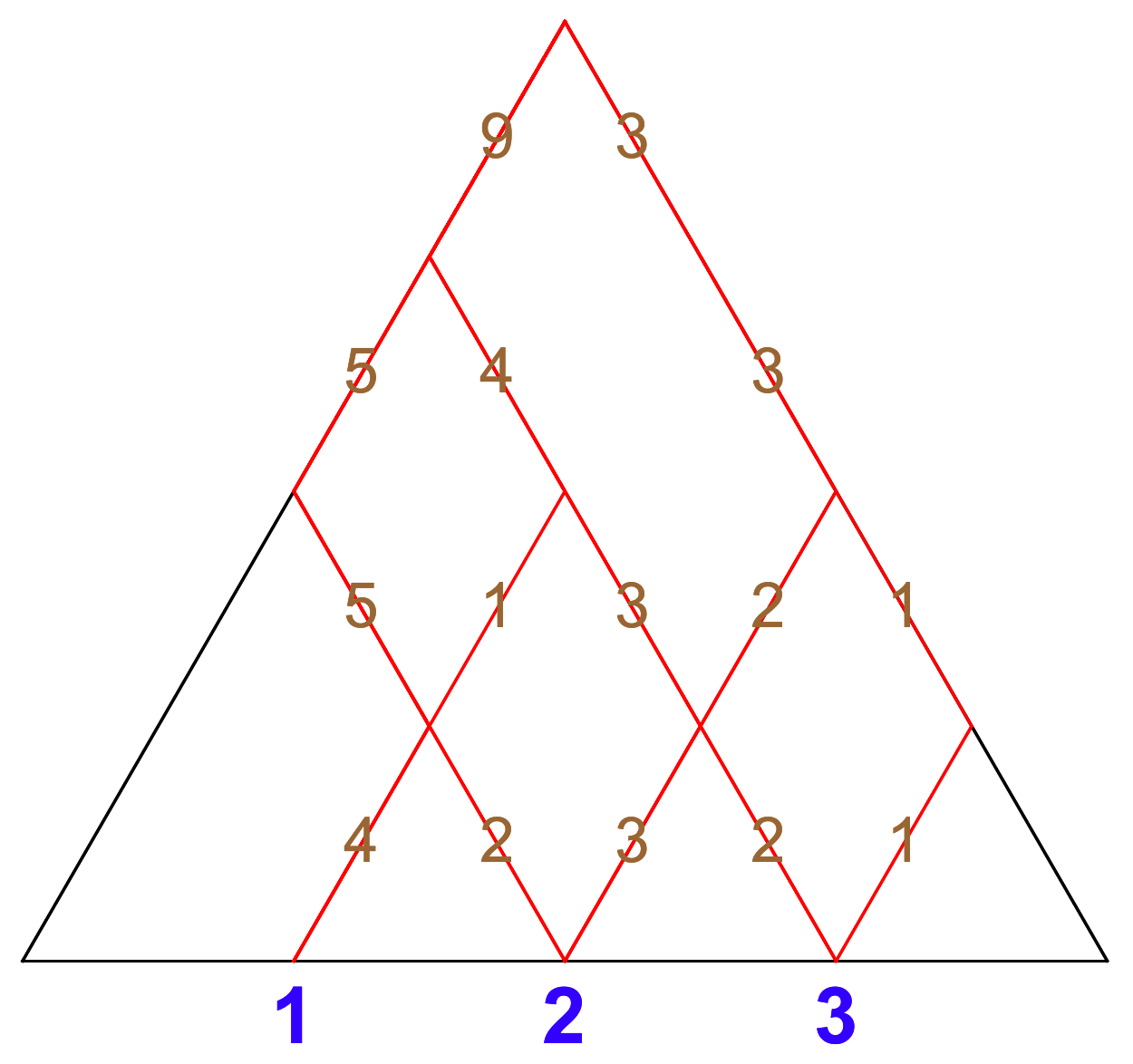}
\caption{SU(4):  Liana forest associated with the Young tableau chosen in the text.\\
Left: With the origin of lianas marked as a superscript.\\ Right: With the origin removed. The corresponding O-blade is in position 19 of Fig.~\ref{allthereducedObladesNew}.} 
\label{lianaForest}
 \end{figure}

The chosen dominant weight (the highest weight $\lambda=\{4,5,3\}$), with associated Young diagram of shape $(12=4+5+3, 8=5+3, 3,0)$, a (decreasing) partition of $23$,  appears at the first level of the liana forest. 
The weight $\xi$ can be read from the sequence giving the total number of left-going edges, when one moves downward, level after level;  namely $(3, 7, 9=8+1, 4)$,
a (non-decreasing) partition of $23$, for which the associated weight is indeed $\delta = \{-4, -2, 5\}=\{3-7, 7-9, 9-4\}$.   This could also be read from the Young tableau given at the beginning of sect.~\ref{ObladesVersusLianas}.
Notice that the sequence of lianas (given by Fig.~\ref{lianas}) or the {\sl indexed liana forest} (that keeps track of the origin of each liana, see Fig.~\ref{lianaForest}, Left) has the same information contents as the Young tableau we started from. 
Part of this information is lost when we superimpose the lianas and remove the indices that give the point they grow from (Fig.~\ref{lianaForest}, Right): several distinct sequences of lianas (\ie distinct Young tableaux, not necessarily semi-standard,  or distinct indexed liana forests) may give rise to the same (unindexed) liana forest. 
If we are only interested in multiplicities and not in the construction of explicit basis vectors in representation spaces, one may forget about {\sl indexed} liana forests because of the existence of a one-to-one correspondence between reduced O-blades and liana forests (compare for instance Fig.~\ref{ObladeSchurExample}, (Middle) and Fig.~\ref{lianaForest}, (Right)).
The number of liana forests being the same as the number of reduced O-blades (in our example they are all displayed on Fig.~\ref{allthereducedObladesNew}), it is also equal to the dimension of the weight subspace defined by $\delta$ in the representation space of highest weight $\lambda$. 
A more detailed analysis of the combinatorics underlying these constructions clearly falls beyond the scope of the present paper,  and we shall stop here.

\smallskip

This way of encoding Young tableaux was explained to one of us (R.C.), more than ten years ago, by A. Ocneanu \cite{AO_Oblades}, who also invented the ``O-blades'' to display the intertwiners that appear in the combinatorics of LR coefficients.
Other aspects of the above forests are summarized in a video lecture: see \cite{AO: VideoGelfand}.  The terminology ``lianas'' is ours.
It was a pleasure (but not so much of a surprise) to rediscover this particular kind of graphical encoding in our discussion of the Schur problem.


\begin{figure}[htb]
 \centering
 \includegraphics[width=32pc]{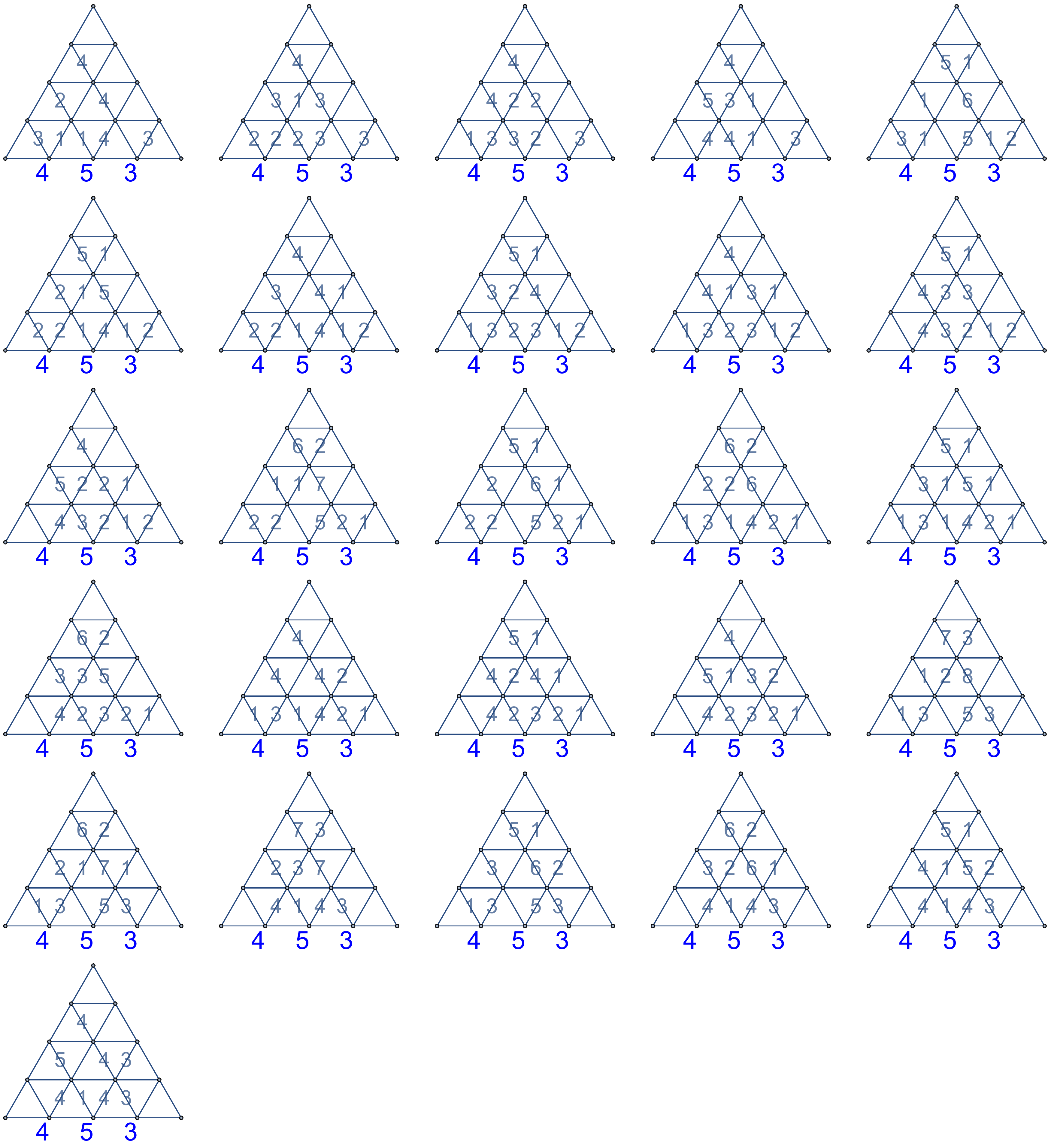}
\caption{SU(4):  All the reduced O-blades for the example discussed in the text} 
\label{allthereducedObladesNew}
 \end{figure}

  \end{document}